\documentclass[a4paper,10pt]{article}
\usepackage{graphicx}
\usepackage{amssymb}

\title{An analytic approach to the Riemann hypothesis}
\date{\today}

\author{Paolo D'Isanto \\
Dipartimento di Fisica ``Ettore Pancini'', \\
Complesso Universitario di Monte S. Angelo, \\ 
Via Cintia, Edificio 6, 80126 Napoli, Italy \and
Giampiero Esposito \\
Dipartimento di Fisica ``Ettore Pancini'', \\
Complesso Universitario di Monte S. Angelo, \\ 
Via Cintia, Edificio 6, 80126 Napoli, Italy \and
INFN Sezione di Napoli, \\
Complesso Universitario di Monte S. Angelo, \\
Via Cintia, Edificio 6, 80126 Napoli, Italy}

\begin{document}
\maketitle
\begin{abstract}
In this work we consider an equation for the Riemann $\zeta$-function 
in the critical half-strip
$$
{\mathbb S}^{+} = \left\{ 
   s = x + iy \in {\mathbb C} \; : \; 
   {\frac{1}{2} < x < 1},  \; 
   {y > 0}  \right \} .
$$ 
With the help of this equation we prove that finding non-trivial zeros of the Riemann
$\zeta$-function outside the critical line ${\rm Re}(s)={1 \over 2}$ would be equivalent to 
the existence of complex numbers $s=x+iy \in {\mathbb S}^{+}$ for which
$$
\sum_{k=1}^{\infty}(-1)^{k+1}{1 \over \zeta(2x)} 
\sum_{n=1}^{\infty}
{\cos \left[y \log \left(1+{k \over n}\right)\right]\over 
n^{2x} \left(1+{k \over n}\right)^{x}}={1 \over 2}.
$$
Such a condition is studied, and the attempt of proving the Riemann hypothesis is found to
involve also the functional equation
$$
\chi_{n}(t)=-\chi_{n} \left(t+{1 \over n} \right),
$$
where $t$ is a real variable $\geq 1$, and $n$ is any natural number.  
The limiting behaviour of the solutions $\chi_{n}(t)$ as $t$ approaches $1$ is then
studied in detail.
\end{abstract}

\section{Introduction}
\setcounter{equation}{0}

This introduction, being written for the general reader, describes the early
work from Euler's definition to Riemann's article.
Riemann's $\zeta$-function is the analytic extension to the whole complex plane of the $\zeta$-function 
\begin{equation}
{\zeta \left( s \right) = \sum_{n = 1}^\infty  {\frac{1}{{n^s  }}} }, 
\label{(1.1)}
\end{equation}
defined by Euler in the region of the complex plane
\begin{equation}
{\mathbb A} = \left\{ {s \in {\mathbb C}:{\rm Re}(s) > 1} \right\} ,
\label{(1.2)}
\end{equation}
where it is absolutely convergent\footnote{This region is an IP-set since if 
$s_1,s_2\in{\mathbb A}$ then $s_1+s_2\in {\mathbb A}$. This observation will be very important in 
the next sections, where we will prove a fundamental identity for  
$\zeta(s_{1}+s_{2})$.}. 
Euler proved that in the region ${\mathbb A}$ the 
$\zeta$-function admits a product representation
\begin{equation}
\zeta (s) = \prod_{p\in {\mathbb P}}  
{\frac{1}{{\left(1 - \frac{1}{{p^s }}\right)}}} \; {\rm with} \; s \in {\mathbb A},  \\
\label{(1.3)}
\end{equation}
where ${\mathbb P}$ is the set of prime numbers\footnote{For a rigorous proof of this relation,
see for example the book by Schwartz \cite{rif1}.}. 
By virtue of the previous relation Euler was 
able to prove that in the region ${\mathbb A}$ the $\zeta$-function has neither zeros nor poles. Riemann 
was inspired by this observation to write his masterpiece article\footnote{Riemann, On the number 
of primes less than a given quantity \cite{rif2}.} in which he looked for an analytical expression 
for the step-function $\pi\left(x\right)$ which counts the number of primes less than a given 
number $x\in {\mathbb R}^{+}$. During his investigation, he discovered three very important 
properties of the $\zeta$-function, probably noted but not proved before by Euler. 
The first was that the $\zeta$-function can be analytically extended to the whole complex 
plane by virtue of the integral representation

\begin{equation}
\zeta \left( s \right) = \frac{{\Gamma \left( {1 - s} \right)}}
{{2\pi i}}\int\limits_{ \gamma} {\frac{{z^{s - 1} }}{{\left(e^z  - 1\right)}}dz}, 
\label{(1.4)}
\end{equation}

where the symbol $\gamma$ indicates a Hankel's contour of the kind in Fig. 1.

\begin{figure}
\centering
\includegraphics{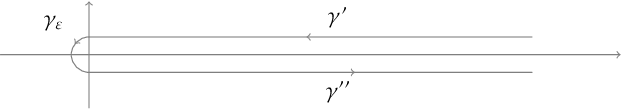}
\caption{Hankel contour for the integral (1.4).}
\end{figure}

The second was the discovery of the functional equation\footnote{In 
the literature, it is written also in the form \cite{rif3}
\[
\zeta \left( s \right) = \Gamma \left( {1 - s} \right)\left( {2\pi } \right)^{s - 1} 
\sin \left( {\frac{{\pi s}}{2}} \right)\zeta \left( {1 - s} \right)
\]}
\begin{equation}
\pi^{-{s \over 2}} \Gamma \left(s \over 2\right)\zeta\left(s\right)
=\pi^{-\frac{1-s}{2}} \Gamma \left(1-s\over 2\right)\zeta\left(1-s\right),
\label{(1.5)}
\end{equation}
which connects the properties of $\zeta$ in the two half-planes in which the complex plane 
remains divided by the vertical line 
\begin{equation}
r_{1/2} = \left\{s \in {\mathbb C}:{\rm Re}(s)=\frac12\right\},
\label{(1.6)}
\end{equation}
named {\it the critical line}. Riemann 
pointed out that, by virtue of the functional equation (1.5), upon
setting $s_k=-2k$ one has $\zeta\left(s_k\right)=0$ for any $k\in {\mathbb N}$. These zeros, being 
obtained from the $\sin$ function, are called 
trivial zeros. Nevertheless, this relation does 
not assure us that they are all the zeros of the $\zeta$, but it proves only that, if other sets 
of zeros exist, they must lie in the critical strip
\begin{equation}
{\mathbb S} = \left\{ {s = x+iy \in {\mathbb C}:0 < x < 1} \right\},
\label{(1.7)}
\end{equation}
and that the zeros have to be located symmetrically about the critical line. These zeros are 
called {\it non-trivial zeros}.\footnote{In the next sections we divide the critical 
strip ${\mathbb S}$ into four parts: two above the real axis
\begin{equation}
\begin{array}{*{20}c}
   {{}^+{\mathbb S}   = \left\{ {s \in {\mathbb C}:{\begin{array}{*{20}c}
   {0 < {\rm Re}(s) < \frac{1}{2}}  \\
   {{\rm Im}(s) > 0}  \\
\end{array}}} \right\}} 
& {\rm and} & {{\mathbb S}^+   = \left\{ {s \in {\mathbb C}: {\begin{array}{*{20}c}
   {\frac12 < {\rm Re}(s) < 1}  \\
   {{\rm Im}(s) > 0}  \\
\end{array}} } \right\}}  \\
\end{array}
\label{(1.8)}
\end{equation}
and two below it
\begin{equation}
\begin{array}{*{20}c}
   {{}_-{\mathbb S}   = \left\{ {s \in {\mathbb C}: {\begin{array}{*{20}c}
   {0 < {\rm Re}(s) < \frac12}  \\
   {{\rm Im}(s) < 0}  \\
\end{array}} } \right\}} 
& {\rm and} & {{\mathbb S}_-   = \left\{ {s \in {\mathbb C}: {\begin{array}{*{20}c}
   {\frac12 < {\rm Re}(s) < 1}  \\
   {{\rm Im}(s) < 0}  \\
\end{array}} } \right\}}  \\
\end{array}
\label{(1.9)}
\end{equation}
and we will refer to ${\mathbb S}^+$ as the {\it critical half-strip} because, if a zero belongs to 
${\mathbb S}^+$, the functional equation imposes the existence of a twin zero in 
${\mathbb S}_-$, and the same holds in $^+{\mathbb S}$ and $_-{\mathbb S}$. 
Now, for any $\varepsilon>0$, we introduce the $\varepsilon$-contraction of the critical 
half-strip ${\mathbb S}_\varepsilon^+$ and the compact $\varepsilon$-contraction of the critical 
half-strip $\left[{\mathbb S}^{\pm}_{\varepsilon,T}\right]$ below a given quantity as
\begin{eqnarray}
\; & \; &
{\mathbb S}_\varepsilon^{+}   = \left \{ 
s \in {\mathbb C}:
\frac{1}{2} + \varepsilon  < {\rm Re}(s) < 1 - \varepsilon, \;
{\rm Im}(s) > 0  \right \} ,
\nonumber \\ 
& \; & 
[{\mathbb S}_{\varepsilon,T}^{+}] 
= \left \{ 
s \in {\mathbb C}:
\frac{1}{2} + \varepsilon  \leq {\rm Re}(s) \leq 1 - \varepsilon, \; 
0 \leq {\rm Im}(s) \leq T \right \}.  
\label{(1.10)}
\end{eqnarray}}

\begin{figure}
\centering
\includegraphics{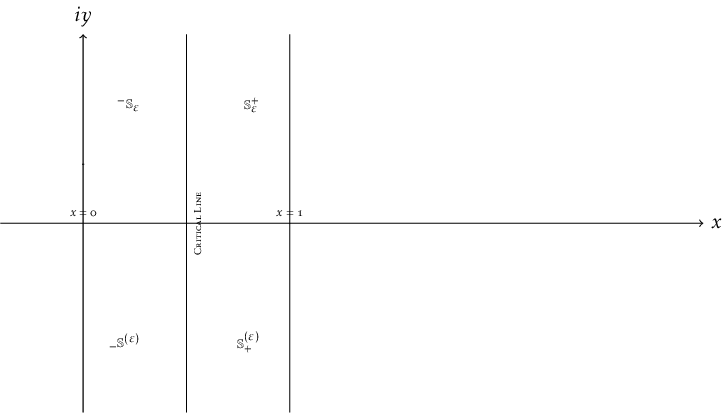}
\caption{The four parts of the critical strip.}
\end{figure}

As far as their existence is concerned, Hardy \cite{HAR1914} proved that, 
on the Critical Line $r_{1/2}$, there exist infinitely many zeros, but
\begin{quote}
The Riemann Hypothesis states that there are no non-trivial zeros outside the critical line $r_{1/2}$.
\end{quote}
The third result obtained by Riemann was the discovery of a \textsc{Product Representation} 
in $[{\mathbb S}_\varepsilon]$ of the function
\begin{equation}
\xi (s) = \frac{{s\left( {s - 1} \right)}}{2}\pi ^{ - \frac{s}{2}} \Gamma 
\left( {\frac{s}{2}} \right)\zeta \left( s \right),
\label{(1.11)}
\end{equation}
symmetric, about the Critical Line, of the form
\begin{equation}
\xi (s) = \xi \left( 0 \right)\prod_{\rho  \in {\cal R}} 
{\left( {1 - \frac{s}{\rho }} \right)} ,
\label{(1.12)}
\end{equation}
where ${\cal R}$ is the set of all non-trivial zeros of the $\zeta$-function. By virtue of 
this representation the Functional Equation (1.5) becomes $\xi\left(s\right)=\xi\left(1-s\right)$. 
But how are they related with prime numbers? As we said before, the goal of Riemann's masterpiece article 
was an analytic expression for the step-function $\pi\left(x\right)$. After having discovered these 
three properties he proposes the definition of the step function
$$
J(x) = \sum_{n=1}^{\infty}  
{\frac{1}{n}\pi \left( {x^{1/n} } \right)}, 
$$
which is, for any $x\in {\mathbb R}^{+}$, a finite sum, since by definition $\pi\left(\alpha\right)\equiv 0$ 
for any $\alpha<2$. Thus, setting $x^{1/n}<2$ we have, for any $n>\frac{\log x}{\log 2}$, that 
$\pi\left(x^{1/n}\right)\equiv 0$. This sum can be inverted by the \textsc{M\"{o}bius Inversion Formula}, 
obtaining\footnote{The M\"{o}bius $\mu$-function is defined, when the number $n$ 
is expressed in the form
\[
n = \prod\limits_{k = 1}^\infty  {p_k^{\alpha _{nk} } }, 
\]
as the function that vanishes when one of the $\alpha_k$ is bigger than one:
\[
\mu(n) = \left\{ {\begin{array}{*{20}c}
   1 & {\rm iff} & {\alpha _{n, i_1}  =  \cdots  = \alpha _{n,i_{2k}}  = 1}  \\
   0 & {} & {\exists \alpha _{n, i_l}  > 1}  \\
   { - 1} & {} & {\alpha _{n, i_1}  =  \cdots  = \alpha _{n,i_{2k + 1}}  = 1}  \\
\end{array}} \right.
\]}
\begin{equation}
\pi \left( x \right) = \sum_{n=1}^{\left[ {\frac{{\log x}}{{\log 2}}} \right]} 
{\frac{{\mu \left( n \right)}}{n}J\left( {x^{1/n} } \right)}. 
\label{(1.13)}
\end{equation}
His main result follows by a careful use of (1.3), hence obtaining the analytic 
expression of $J\left(x\right)$ as 
(on denoting by ${\rm Li}(x)= \int_{2}^{x} {dt \over (\log t)}$ the logarithmic integral)
\begin{equation}
J\left( x \right) = {\rm Li}(x) - \sum\limits_{\rho  \in {\cal R}} {{\rm Li}\left( {x^\rho  } \right)}  
- \log 2 + \int\limits_x^\infty  {\frac{{dt}}{{t\left( {t^2  - 1} \right)\log t}}}. 
\label{(1.14)}
\end{equation}
Upon substituting this expression into (1.13) he obtained the desired result. The dominant term in the 
\textsc{Riemann's main formula} (1.14) is
\[
{\rm Li}\left(x\right)  
= \frac{x}{{\log x}} -{2 \over \log 2}
+\int_{2}^{x} {dt \over (\log t)^{2}},
\]
which is the asymptotic estimate of the number of primes between 2 and $x$, as stated by the 
\textsc{Prime Number Theorem} (here $a$ is a positive constant):
\begin{equation}
\left | \pi(x)-{\rm Li}(x) \right | = {\rm O}\left(x e^{-a \sqrt{\log(x)}} \right),
\label{(1.15)}
\end{equation}
proved by Hadamard and de La Vall\'ee-Poussin \cite{PNT1,PNT2}. 
It was later proved by von Koch \cite{VK} that the remainder term in (1.15) is 
${\rm O}(\sqrt{x}\log(x))$ if and only if the Riemann hypothesis holds.
But the most interesting term of (1.14) 
is the second, whose sum runs over $\rho \in {\cal R}$. Thus, Riemann understood that his dream to obtain 
a staircase-function, whose steps were localized on the prime numbers, could be realized if and only 
if the positions of non-trivial zeros of the Riemann's $\zeta$-function were known. 

During the twentieth century, in the forties, it became
clear that the Riemann $\zeta$-function is an element of a larger class of 
functions, called $L$-functions \cite{ANT}, but in our work we deal with
Riemann's $\zeta$-function only. Our paper is entirely devoted to a {\it purely analytic
approach to the Riemann hypothesis}. For this purpose,
Sec. $2$ outlines the strategy we adopt; Sec. $3$ studies 
a double series in the critical half-strip ${\mathbb S}^{+}$ that is crucial for our investigation;
Sec. $4$ obtains an equation satisfied by Riemann's $\zeta$-function in such a critical
half-strip; Sec. $5$ obtains a necessary condition for finding non-trivial zeros in the
critical half-strip ${\mathbb S}^{+}$; Sec. $6$ exploits the result in Sec. $5$ 
and derives the conditions under which
the existence of non-trivial zeros in ${\mathbb S}^{+}$ leads to a contradiction; an 
assessment of our approach is presented in Sec. $7$, and some important technical results
are provided by the $3$ Appendices. 

\section{The strategy we adopt}
\setcounter{equation}{0}

In order to prove that there are no non-trivial zeros outside the critical line,
we consider an analytic continuation of the Riemann $\zeta$-function (1.1) 
to the whole critical strip ${\mathbb S}$ by means of the series with alternating signs
\begin{equation}
\zeta(s) = {1 \over (1-2^{1-s})} \sum_{n=1}^{\infty}
{(-1)^{n-1}\over n^{s}}, \; s \in {\mathbb A} \cup {\mathbb S}.
\label{(2.1)}
\end{equation}
On the right-hand side of Eq. (2.1) we recognize the Dirichlet $\eta$-function
\begin{equation}
\eta(s) = 1-{1 \over 2^{s}}+{1 \over 3^{s}} - ...
=\sum_{n=1}^{\infty}{(-1)^{(n-1)}\over n^{s}},
\label{(2.2)}
\end{equation}
and one can easily prove that its zeros within the open strip, i.e. the critical strip
${\mathbb S}$ deprived of the vertical lines $s=0$ and $s=1$, coincide with non-trivial
zeros of the Riemann $\zeta$-function. If we can now find, for all $\varepsilon >0$, a
{\it necessary} condition for the existence of zeros in ${\mathbb S}_{\varepsilon}^{+}$
(i.e. the half-strip on the top right-hand sector of ${\mathbb S}$), for the 
$\eta$-function the invalidation of such a condition will be {\it sufficient} to prove
the lack of zeros in ${\mathbb S}_{\varepsilon}^{+}$ for the $\zeta$-function, and hence
the non-existence of zeros for every analytic continuation of the $\zeta$-function.
In order to obtain the desired necessary condition mentioned above, we obtain
preliminarily a fundamental identity by introducing a convergent double series in 
${\mathbb S}^{+}$, where such a series is only conditionally convergent. 

\section{A double series in ${\mathbb S}^+$}
\setcounter{equation}{0}

The first step is to study the double series\footnote{In the Appendix B we give the proof of 
the Pringsheim-convergence (a concept defined in Eq. (A.6)) of the next series.}
\begin{equation}
\sum\limits_{n_1  \ne n_2 } {\frac{{\left( { - 1} \right)^{n_1  + n_2 } }}{{n_1^s n_2^{\bar s} }}}, 
\label{(3.1)}
\end{equation}
where $n_1,n_2 \in {\mathbb N}$, and $s \equiv x+iy$, $\bar s \equiv x-iy$ are defined in 
${\mathbb S}^+$ where $2x>1$. Hence we can write
\begin{eqnarray}
\; & \; &
\sum\limits_{n_1  \ne n_2 } {\frac{{\left( { - 1} \right)^{n_1  + n_2 } }}{{n_1^s n_2^{\bar s} }}}  
= \sum\limits_{n_1  \ne n_2 } {\frac{{\left( { - 1} \right)^{n_1  + n_2 } }}{{\left( {n_1 n_2 } 
\right)^x }}} \left( {\frac{{n_2 }}{{n_1 }}} \right)^{iy} 
\nonumber \\
& = & \sum\limits_{n_1  \ne n_2 } 
{\frac{{\left( { - 1} \right)^{n_1  + n_2 } }}{{\left( {n_1 n_2 } \right)^x }}
e^{iy\log \frac{{n_2 }}{{n_1 }}} }   
\nonumber \\
&=& \sum\limits_{n_1  \ne n_2 } {\frac{{\left( { - 1} \right)^{n_1  + n_2 } }}
{{\left( {n_1 n_2 } \right)^x }}\left\{ {\cos \left( {y\log \frac{{n_2 }}{{n_1 }}} \right) 
+ i\sin \left( {y\log \frac{{n_2 }}{{n_1 }}} \right)} \right\}}   
\nonumber \\
&=& \sum\limits_{n_1  \ne n_2 } \frac{{\left( { - 1} \right)^{n_1  + n_2 } }}{{\left( {n_1 n_2 } 
\right)^x }}\cos \left( {y\log \frac{{n_2 }}{{n_1 }}} \right) 
\nonumber \\
& + & i\sum\limits_{n_1  \ne n_2 } 
{\frac{{\left( { - 1} \right)^{n_1  + n_2 } }}{{\left( {n_1 n_2 } \right)^x }}\sin 
\left( {y\log \frac{{n_2 }}{{n_1 }}} \right)} .  
\label{(3.2)}
\end{eqnarray}
Since the sum over $n_{1} \not = n_{2}$ can be split into $2$ sums according to
\begin{equation}
\sum\limits_{n_1  \ne n_2 } {\left(  \cdot  \right)}  = \sum\limits_{n_1  > n_2 } 
{\left(  \cdot  \right)}  + \sum\limits_{n_2  > n_1 } {\left(  \cdot  \right)}, 
\label{(3.3)}
\end{equation}
we can rewrite (3.2) as
\begin{eqnarray}
\; & \; &
\sum\limits_{n_1  \ne n_2 } {\frac{{\left( { - 1} \right)^{n_1  + n_2 } }}
{{n_1^s n_2^{\bar s} }}}
\nonumber \\   
&=& \biggr[ \sum\limits_{n_1  > n_2 } {\frac{{\left( { - 1} \right)^{n_1  + n_2 } }}
{{\left( {n_1 n_2 } \right)^x }}\cos \left( {y\log \frac{{n_2 }}{{n_1 }}} \right)} 
\nonumber \\
& + & \sum\limits_{n_2  > n_1 } {\frac{{\left( { - 1} \right)^{n_1  + n_2 } }}{{\left( {n_1 n_2 } 
\right)^x }}\cos \left( {y\log \frac{{n_2 }}{{n_1 }}} \right)} \biggr]   
\nonumber \\
&+& i \biggr [ \sum\limits_{n_1  > n_2 } {\frac{{\left( { - 1} \right)^{n_1  + n_2 } }}
{{\left( {n_1 n_2 } \right)^x }}\sin \left( {y\log \frac{{n_2 }}{{n_1 }}} \right)}
\nonumber \\
&+& \sum\limits_{n_2  > n_1 } {\frac{{\left( { - 1} \right)^{n_1  + n_2 } }}{{\left( {n_1 n_2 } 
\right)^x }}\sin \left( {y\log \frac{{n_2 }}{{n_1 }}} \right)}  \biggr ].
\label{(3.4)}
\end{eqnarray}
Since in the exchange $n_1\leftrightarrow n_2$ we have
\begin{eqnarray}
\; & \; &
{\sum\limits_{n_2  > n_1 } {\frac{{\left( { - 1} \right)^{n_1  + n_2 } }}{{\left( {n_1 n_2 } 
\right)^x }}\cos \left( {y\log \frac{{n_2 }}{{n_1 }}} \right)}  
= \sum\limits_{n_1  > n_2 } 
{\frac{{\left( { - 1} \right)^{n_1  + n_2 } }}{{\left( {n_1 n_2 } \right)^x }}\cos 
\left( {y\log \frac{{n_1 }}{{n_2 }}} \right)}} , 
\nonumber \\
& {} &   
\sum\limits_{n_2  > n_1 } {\frac{{\left( { - 1} \right)^{n_1  + n_2 } }}{{\left( {n_1 n_2 } 
\right)^x }}\sin \left( {y\log \frac{{n_2 }}{{n_1 }}} \right)}  
\nonumber \\
& = & \sum\limits_{n_1  > n_2 } 
{\frac{{\left( { - 1} \right)^{n_1  + n_2 } }}{{\left( {n_1 n_2 } \right)^x }}
\sin \left( {y\log \frac{{n_1 }}{{n_2 }}} \right)}  
\label{(3.5)}
\end{eqnarray}
together with the obvious symmetries
\begin{equation}
{\cos \left( {y\ln \frac{{n_2 }}{{n_1 }}} \right) 
= \cos \left( {y\ln \frac{{n_1 }}{{n_2 }}} \right)}, \;  
{\sin \left( {y\ln \frac{{n_2 }}{{n_1 }}} \right) 
=- \sin \left( {y\ln \frac{{n_1 }}{{n_2 }}} \right)},  
\label{(3.6)}
\end{equation}
we can write
\begin{eqnarray}
\; & \; & \sum\limits_{n_1  \ne n_2 } {\frac{{\left( { - 1} \right)^{n_1  + n_2 } }}
{{n_1^s n_2^{\bar s} }}}  \nonumber \\ 
&=& \left\{ {\sum\limits_{n_1  > n_2 } {\frac{{\left( { - 1} \right)^{n_1  + n_2 } }}
{{\left( {n_1 n_2 } \right)^x }}\cos \left( {y\log \frac{{n_2 }}{{n_1 }}} \right)}  
+ \sum\limits_{n_1  > n_2 } {\frac{{\left( { - 1} \right)^{n_1  + n_2 } }}
{{\left( {n_1 n_2 } \right)^x }}\cos \left( {y\log \frac{{n_1 }}{{n_2 }}} \right)} } \right\}
\nonumber \\ 
&+& i\left\{ {\sum\limits_{n_1  > n_2 } {\frac{{\left( { - 1} \right)^{n_{1}+n_{2}} }}
{{\left( {n_1 n_2 } \right)^x }}\sin \left( {y\log \frac{{n_2 }}{{n_1 }}} \right)}  
+ \sum\limits_{n_1  > n_2 } {\frac{{\left( { - 1} \right)^{n_{1}+n_{2}} }}
{{\left( {n_1 n_2 } \right)^x }}
\sin \left( {y\log \frac{{n_1 }}{{n_2 }}} \right)} } \right\}  
\nonumber \\
&=& \sum\limits_{n_1  > n_2 } {\frac{{\left( { - 1} \right)^{n_1  + n_2 } }}
{{\left( {n_1 n_2 } \right)^x }}} \left\{ {\cos \left( {y\log \frac{{n_2 }}{{n_1 }}} \right) 
+ \cos \left( {y\log \frac{{n_1 }}{{n_2 }}} \right)} \right\} 
\nonumber \\
&+& i\sum\limits_{n_1  > n_2 } 
{\frac{{\left( { - 1} \right)^{n_1  + n_2 } }}{{\left( {n_1 n_2 } \right)^x }}} 
\left\{ {\sin \left( {y\log \frac{{n_2 }}{{n_1 }}} \right) + \sin \left( 
{y\log \frac{{n_1 }}{{n_2 }}} \right)} \right\}  
\nonumber \\
&=& \sum\limits_{n_1  > n_2 } {\frac{{\left( { - 1} \right)^{n_1  + n_2 } }}
{{\left( {n_1 n_2 } \right)^x }}} \left\{ {\cos \left( {y\log \frac{{n_2 }}{{n_1 }}} \right) 
+ \cos \left( {y\log \frac{{n_2 }}{{n_1 }}} \right)} \right\}
\nonumber \\ 
&+& i\sum\limits_{n_1  > n_2 } 
{\frac{{\left( { - 1} \right)^{n_1  + n_2 } }}{{\left( {n_1 n_2 } \right)^x }}} 
\left\{ {\sin \left( {y\log \frac{{n_2 }}{{n_1 }}} \right) - \sin 
\left( {y\log \frac{{n_2 }}{{n_1 }}} \right)} \right\}  
\nonumber \\
&=& 2\sum\limits_{n_1  > n_2 } {\frac{{\left( { - 1} \right)^{n_1  + n_2 } }}
{{\left( {n_1 n_2 } \right)^x }}} \cos \left( {y\log \frac{{n_2 }}{{n_1 }}} \right).  
\label{(3.7)}
\end{eqnarray}
For any $n,k \in {\mathbb N}$, if we define $n_1\equiv n+k$ and $n_2\equiv n$ we can 
rewrite the previous relation in the form
\begin{eqnarray}
\; & \; &
\sum\limits_{n_1  \ne n_2 } {\frac{{\left( { - 1} \right)^{n_1  + n_2 } }}
{{n_1^s n_2^{\bar s} }}}  
=2\sum\limits_{n + k > n} {\frac{{\left( { - 1} \right)^{n + k + n} }}{{\left[ 
{n\left( {n + k} \right)} \right]^x }}\cos \left[ {y\log \left( {\frac{{n + k}}{n}} 
\right)} \right]}  
\nonumber \\
&=& 2\sum\limits_{k = 1}^\infty  {\left( { - 1} \right)^k 
\sum\limits_{n = 1}^\infty  {\frac{{\cos \left[ {y\log \left( {1 + \frac{k}{n}} \right)} 
\right]}}{{\left[ {n\left( {n + k} \right)} \right]^x }}} }.   
\label{(3.8)}
\end{eqnarray}
Now, $\zeta\left(2x\right)$ being bounded from above in ${\mathbb S}^+$ we can define the 
{\it zeros' functions}
\begin{equation}
Z_{k+1} (x,y) =  \frac{1}{{\zeta \left( {2x} \right)}}
\sum\limits_{n = 1}^\infty  {\frac{{\cos \left[ {y\log \left( {1 + \frac{k}{n}} \right)} \right]}}
{{\left[ {n\left( {n + k} \right)} \right]^x }}} ,
\label{(3.9)}
\end{equation}
and hence we can write
\begin{eqnarray}
\sum\limits_{n_1  \ne n_2 } {\frac{{\left( { - 1} \right)^{n_1  + n_2 } }}
{{n_1^s n_2^{\bar s} }}}   
&=&  - 2\zeta \left( {2x} \right)\sum\limits_{k = 1}^\infty  {\left( { - 1} \right)^{k + 1} 
\left\{ {\frac{1}{{\zeta \left( {2x} \right)}}\sum\limits_{n = 1}^\infty  {\frac{{\cos 
\left[ {y\log \left( {1 + \frac{k}{n}} \right)} \right]}}{{\left[ {n\left( {n + k} \right)} 
\right]^x }}} } \right\}}  
\nonumber \\
&=& - 2\zeta \left( {2x} \right)\sum\limits_{k = 1}^\infty  
{\left( { - 1} \right)^{k + 1} Z_{k + 1} \left( {x,y} \right)}.   
\label{(3.10)}
\end{eqnarray}
We refer the reader to Appendix B for a crucial remark concerning the $Z_{k+1}$ functions
that we have just introduced.

\section{Riemann's $\zeta$-function in the critical half-strip 
${\mathbb S}^+$} 
\setcounter{equation}{0}

Before reverting to the critical half-strip ${\mathbb S}^+$ we want to prove 
that an identity for the $\zeta$-function holds on 
${\mathbb A}$ (see (1.2)), whose maximal extension 
to the critical half-strip will be useful in the next sections. We start by taking in 
${\mathbb A}$ two different points $s_1,s_2 \in {\mathbb A}$. 
Of course, $s_1+s_2$ lies in ${\mathbb A}$ as well, as we already pointed out 
in the Introduction. Thus, we can prove a lemma as follows.
\vskip 0.3cm
\noindent
{\bf Lemma 1 (Fundamental Identity)}.
For every $s_1$ and $s_2$ in ${\mathbb A}$, for which the series expressing $\zeta(s)$ 
is absolutely convergent, the relation
\begin{equation}
\zeta (s_1 )\zeta (s_2 ) - \zeta (s_1  + s_2 ) = \sum\limits_{n_1  \ne n_2 } 
{\frac{1}{{n_1^{s_1 } n_2^{s_2 } }}} 
\label{(4.1)}
\end{equation}
holds.
\vskip 0.3cm
\noindent
{\it Proof}. Starting from the definition of $\zeta(s)$ in ${\mathbb A}$, 
any series that we write is absolutely convergent, and in particular
\begin{equation}
\left( {\sum\limits_{n_1} {\frac{1}{{n_1^{s_1 } }}} } \right)\left( 
{\sum\limits_{n_2} {\frac{1}{{n_2^{ s_2 } }}} } 
\right) = \sum\limits_{n=n_1=n_2} {\frac{1}{{n^{\left( {s_1  + s_2 } \right)} }}}  
+ \sum\limits_{n_1  \ne n_2 } {\frac{1}{{n_1^{s_1 } n_2^{ s_2 } }}}, 
\label{(4.2)}
\end{equation}
hence the thesis follows.

At this stage we need to revert to the Critical Half Strip 
${\mathbb S}^+$, looking for an equation
valid perhaps now in ${\mathbb S}^+$. Our aim is to obtain, thanks to such an equation, a relation 
in which the non-trivial zeros (i.e. those $s$ not given by $-2k$ with $k$ a positive
integer, as we discussed after Eq. (1.5)) 
are involved in a characterization formula which may be a 
necessary condition for their existence. Upon violating it, we will obtain a sufficient condition for 
the impossibility to have zeros in ${\mathbb S}^+$. In order to extend the Fundamental Identity 
to ${\mathbb S}^+$ we need the proof of the following lemma:
\vskip 0.3cm 
{\bf Lemma 2 (Maximal Extension of the Fundamental Identity)}.
For every $s$ and $\bar s$ with ${\rm Re}(s)>1/2$ we have the relation
\begin{equation}	
\left( {1 - 2^{1 - s} } \right)\left( {1 - 2^{1 - \bar s} } \right)\zeta (s)\zeta (\bar s) 
- \left. {\zeta (2x)} \right|_{2x > 1}  = \sum\limits_{n_1  \ne n_2 } {\frac{{\left( 
{ - 1} \right)^{n_1  + n_2 } }}{{n_1^s n_2^{\bar s} }}}.  
\label{(4.3)}
\end{equation}
\vskip 0.3cm
\noindent
{\it Proof}. We can extend the Fundamental Identity to the critical half-strip by virtue of the 
$\zeta$'s representation (2.1), valid for ${\rm Re}(s)>0$.
Thus, taking $s=x+iy \in {\mathbb S}^{+}$, we have $2x>1$ and hence
$\zeta\left(2x\right)<\infty$, obtaining therefore
\begin{eqnarray}
\; & \; &
\left( {\sum\limits_{n = 1}^\infty  {\frac{{\left( { - 1} \right)^{n - 1} }}{{n^s }}} } 
\right)\left( {\sum\limits_{n = 1}^\infty  {\frac{{\left( { - 1} \right)^{n - 1} }}
{{n^{\bar s} }}} } \right)
= \sum\limits_{n=1}^\infty  {\frac{{\left( { - 1} \right)^{2(n - 1)} }}{{n^{s + \bar s} }}}  
+ \sum\limits_{n_1  \ne n_2 } {\frac{{\left( { - 1} \right)^{n_1  + n_2 } }}{{n_1^s n_2^{\bar s} }}}
\nonumber \\  
&=& \sum\limits_{n=1}^\infty  {\frac{1}{{n^{2x} }}}  
+ \sum\limits_{n_1  \ne n_2 } {\frac{{\left( { - 1} 
\right)^{n_1  + n_2 } }}{{n_1^s n_2^{\bar s} }}}.
\label{(4.4)}
\end{eqnarray}
By re-expressing the first line of Eq. (4.4) with the help of (2.1), we obtain eventually 
the desired Eq. (4.3)

\section{A necessary condition for the zeros of the Riemann $\zeta$-function in the critical half-strip 
${\mathbb S}^+$}
\setcounter{equation}{0}

In this section we obtain a \textsl{necessary condition} for non-trivial zeros in the critical 
half-strip ${\mathbb S}^{+}$. 
In order to achieve our goal, we need to prove the next theorem.
\vskip 0.3cm
\noindent
{\bf Theorem 1 (Critical Half-Strip's Zeros)}. If $s=x+iy$ is a non-trivial Riemann zero 
in the critical half-strip ${\mathbb S}^+$, then (see (3.9))
\begin{equation}
{{\sum\limits_{k = 1}^\infty  {\left( { - 1} \right)^{k + 1} Z_{k+1}\left(x,y\right) } 
= \frac{1}{2}}}.
\label{(5.1)}
\end{equation}
\vskip 0.3cm
\noindent
{\it Proof}. Remembering the Maximal Extension of the Fundamental Identity to the critical 
half-strip ${\mathbb S}^+$
\begin{equation}
\left( {1 - 2^{1 - s} } \right)\left( {1 - 2^{1 - \bar s} } \right)\zeta (s)\zeta (\bar s) 
- \left. {\zeta (2x)} \right|_{2x > 1}  = \sum\limits_{n_1  \ne n_2 } {\frac{{\left( { - 1} 
\right)^{n_1  + n_2 } }}{{n_1^s n_2^{\bar s} }}},  
\label{(5.2)}
\end{equation}
and setting $\zeta\left(s\right)=0$ in it, we can write
\begin{equation}
- \left. {\zeta (2x)} \right|_{2x > 1}  = \sum\limits_{n_1  \ne n_2 } 
{\frac{{\left( { - 1} \right)^{n_1  + n_2 } }}{{n_1^s n_2^{\bar s} }}}. 
\label{(5.3)}
\end{equation}
By virtue of Eq. (3.10) we obtain the thesis.
\vskip 0.3cm
\noindent
{\bf Corollary 1 (Symmetries of zeros about the $x$-axis)}. If $x+iy \in {\mathbb S}^+$ is a 
point of ${\mathbb S}^{+}$ that satisfies the necessary condition for 
non-trivial zeros, the same holds for its complex conjugate $x-iy \in {\mathbb S}_{-}$, i.e.
\begin{equation}
Z_{k+1}(x,y) = Z_{k+1}(x,-y), \; 
x + iy \in {\mathbb S}^{+}  .  		
\label{(5.4)}
\end{equation}
\vskip 0.3cm
\noindent
{\it Proof}. Since the cosine function is even, we find that (see again (3.9))
\begin{equation}
Z_{k+1}(x,y) = Z_{k+1}(x,-y).
\label{(5.5)}
\end{equation}

\section{Th\-e non-tri\-vi\-al ze\-ros of th\-e Rie\-ma\-nn $\zeta$-func\-ti\-on}
\setcounter{equation}{0}
\vskip 0.3cm
\noindent
{\bf Theorem 2 (Riemann hypothesis)}. Any non-trivial zero of the Riemann 
$\zeta$-function has the form $\rho=\frac12+iy$.
\vskip 0.3cm
\noindent
{\it Proof}. We suppose by contradiction that we can find a non-trivial zero $\bar x+i\bar y$ 
outside the critical line. By virtue of Eq. (5.1) we have necessarily that 
\begin{equation}
\sum\limits_{k = 1}^\infty  {\left( { - 1} \right)^{k + 1} Z_{k+1} 
\left( {\bar x,\bar y} \right)}  = \frac{1}{2}.
\label{(6.1)}
\end{equation}
On the other hand, by defining the \textsl{finite-difference} operators \cite{Jordan,Davies}
\begin{eqnarray}
{\bf \Delta} _h f\left( t \right) & = & f\left( {t + h} \right) - f\left( t \right),  
\nonumber \\
\textbf{E}_h f\left( t \right) & = & \left( {I + {\bf \Delta}_h } 
\right)f\left( t \right) = f\left( {t + h} \right),  
\nonumber \\
\textbf{M}_h f\left( t \right) & = & \frac{1}{2}\left( {I + \textbf{E}_h } \right)f
\left( t \right) = \frac{1}{2}\left\{ {f\left( t \right) + f\left( {t + h} \right)} \right\},  
\label{(6.2)}
\end{eqnarray}
ch\-oo\-si\-ng $h=1/n$, $f\left(t\right) = \frac{\cos\left[y\log\left(t\right)\right]}{t^x}$
(whose domain is ${\mathbb R}^{+}$), and po\-in\-ti\-ng o\-ut th\-at 
$\textbf{E}^{k}_{h}f(t) = f(t+kh)$ w\-e ca\-n se\-t
\begin{equation}
\frac{{\cos \left[ {y\log \left( {1 + \frac{k}{n}} \right)} \right]}}{{\left( {1 + \frac{k}{n}} \right)^x }} 
= f\left( {1 + \frac{k}{n}} \right) \equiv \mathop{\lim}\limits_{t\to 1} f\left( {t + \frac{k}{n}} \right) 
= \mathop{\lim}\limits_{t\to 1}\textbf{E}_{1/n}^k f\left( t \right),
\label{(6.3)}
\end{equation}
hence we can write
\begin{equation}
\sum\limits_{k = 1}^\infty  {\left( { - 1} \right)^{k + 1} \frac{{\cos \left[ {y\log 
\left( {1 + \frac{k}{n}} \right)} \right]}}{{\left( {1 + \frac{k}{n}} \right)^x }}}  
= \sum_{k=1}^{\infty} (-1)^{k+1} 
\mathop{\lim}\limits_{t\to 1} f
\left( {t + \frac{k}{n}} \right).
\label{(6.4)}
\end{equation}
At this stage, it is of crucial importance to understand whether, on the right-hand side 
of (6.4), we can bring the limit outside the summation $\sum_{k=1}^{\infty}$. For this
purpose, we need first a definition, and a theorem proved 
hereafter in three steps. Our approach
aims at invalidating the condition (5.1), which proves in turn the Riemann hypothesis.
 
\subsection{Legitimacy of exchanging limit as $t \rightarrow 1$ with summation over $k$} 

\noindent
{\bf Definition}. A sequence of equicontinuous and uniformly bounded functions $S_{L,N}(t)$ 
in the closed interval $[1,A]$ is here said to be of Cauchy-Ces\`aro type with respect to the
subscript $L$ if and only if, for all $\varepsilon >0$, the exists a $\nu_{\varepsilon} \in
{\mathbb N}$ such that, for all $M > \nu_{\varepsilon}$, one has
\begin{equation} 
\left | S_{L+1,N}(t)-S_{L,N}(t) \right |_{(C,1)} =
\left | {1 \over M} \sum_{L=1}^{M} \Bigr(S_{L+1,N}(t)-S_{L,N}(t)\Bigr) \right |
< \varepsilon,
\label{(6.5)}
\end{equation}
for all $t$ in the closed interval $[1,A]$.

By relying upon this definition, our analysis goes on as follows.
\vskip 0.3cm
\noindent
{\bf Theorem 2.0}. The double sequence of equicontinuous and uniformly bounded functions
\begin{equation}
S_{L,N}(t) = \sum_{k=1}^{L}(-1)^{k+1} {1 \over \zeta(2x)}
\sum_{n=1}^{N}{1 \over n^{2x}}
{\cos \left[y \log \left(t+{k \over n}\right)\right] \over 
\left(t+{k \over n}\right)^{x}},
\label{(6.6)}
\end{equation}
defined for $t \in [1,A]$, with $A < \infty$, converges uniformly to the function
\begin{equation}
S(t) = \lim_{L \to \infty} \; \lim_{N \to \infty} S_{L,N}(t)
=\lim_{N \to \infty} \; \lim_{L \to \infty} S_{L,N}(t).
\label{(6.7)}
\end{equation}
\vskip 0.3cm
\noindent
{\it Proof}. We must prove the following majorization:
\begin{eqnarray}
\; & \; & \left | S_{L+1,N+1}(t)-S_{L,N}(t) \right | 
\nonumber \\
& \leq & \left | S_{L+1,N+1}(t)-S_{L+1,N}(t) \right | 
+\left | S_{L+1,N}(t)-S_{L,N+1}(t) \right |
\nonumber \\ 
&+& \left | S_{L,N+1}(t)-S_{L,N}(t) \right | < 3 \varepsilon .
\label{(6.8)}
\end{eqnarray}
{\bf First step}. Upon bearing in mind the definition (6.6) we can write, for the
first term on the second line of (6.8),
\begin{eqnarray}
\; & \; & \Bigr | S_{L+1,N+1}(t)-S_{L+1,N}(t) \Bigr | 
\nonumber \\
&=& \biggr | \sum_{k=1}^{L+1} (-1)^{k+1} \left \{ {1 \over \zeta(2x)}
\sum_{n=1}^{N+1}{1 \over n^{2x}}
{ \cos \left [ y \log \left(t+{k \over n}\right)\right]
\over \left(t+{k \over n}\right)^{x}} \right \}
\nonumber \\
&-& \sum_{k=1}^{L+1} (-1)^{k+1} \left \{
{1 \over \zeta(2x)} \sum_{n=1}^{N} {1 \over n^{2x}}
{\cos \left[y \log \left(t+{k \over n}\right)\right] \over
\left(t+{k \over n}\right)^{x}} \right \} \biggr |
\nonumber \\
&=& {1 \over \zeta(2x)} \biggr |
\sum_{k=1}^{L+1} (-1)^{k+1} \left \{
\sum_{n=1}^{N+1}{1 \over n^{2x}}
{\cos \left[y \log \left(t+{k \over n}\right) \right] \over
\left (t+{k \over n}\right)^{x}} \right . 
\nonumber \\
&-& \left . \sum_{n=1}^{N} {1 \over n^{2x}}
{\cos \left[y \log \left(t+{k \over n}\right)\right]
\over \left(t+{k \over n}\right)^{x}} \right \} \biggr |
\nonumber \\
&=& {1 \over \zeta(2x)} \left | \sum_{k=1}^{L+1} (-1)^{k+1}
{1 \over (N+1)^{2x}}
{\cos \left[y \log \left(t+{k \over (N+1)}\right)\right]
\over \left(t+{k \over (N+1)}\right)^{x}} \right |
\nonumber \\
& \leq & {1 \over \zeta(2x)} \left | \sum_{k=1}^{L+1} (-1)^{k+1} b_{k}^{(N)}(t)
\right | ,
\label{(6.9)}
\end{eqnarray}
where we have defined
\begin{equation}
b_{k}^{(N)}(t) = {1 \over (N+1)^{2x}}
{\cos \left[y \log \left(t+{k \over (N+1)}\right)\right]
\over \left(t+{k \over (N+1)}\right)^{x}}.
\label{(6.10)}
\end{equation}
Thus, for the sequence $S_{L,N}$ to be of Cauchy type with respect to $N$, it is sufficient to
prove that, for all $\varepsilon >0$, there exists a $\mu_{\varepsilon} \in \mathbb{N}$
such that, for all $N > \mu_{\varepsilon}$, one has
\begin{equation}
\sigma = \left | \sum_{k=1}^{L+1} (-1)^{k+1} b_{k}^{(N)}(t) \right | 
< \; \varepsilon .
\label{(6.11)}
\end{equation}
Indeed, on using the standard notation according to which $[a]$ is the integer part of 
a rational number $a$, one finds
\begin{eqnarray}
\; & \; & \sigma = {1 \over (N+1)^{2x}} 
\left | \sum_{k=1}^{[{{L+1} \over 2}]}
{\cos \left[y \log \left(t+{(2k-1) \over N+1)}\right)\right]
\over \left(t+{(2k-1)\over (N+1)}\right)^{x}} \right .
\nonumber \\
& - & \left . \sum_{k=1}^{[{{L+1} \over 2}]}
{\cos \left[y \log \left(t+{2k \over (N+1)}\right)\right]
\over \left(t+{2k \over (N+1)}\right)^{x}} \right | 
\nonumber \\
& \leq & {1 \over (N+1)^{2x}} \sum_{k=1}^{[{{L+1} \over 2}]}
\left | {\cos \left[y \log \left(t+{(2k-1)\over (N+1)}\right)\right]
\over \left(t+{(2k-1)\over (N+1)}\right)^{x}} \right .
\nonumber \\
&-& \left . {\cos \left[y \log \left(t+{2k \over (N+1)}\right)\right]
\over \left(t+{2k \over (N+1)}\right)^{x}} \right |
\label{(6.12)}
\end{eqnarray}
and this shows that it is sufficient to prove that, as $L$ approaches $\infty$,
the resulting series on the third line of (6.12) is convergent for all 
$x+iy \in {\mathbb S}^{+}$. As a matter of fact, by applying the Abel criterion
\cite{Giusti} to the strip 
${\mathbb S}^{+}$, where $x \in ]{1 \over 2},1 [$, we find
\begin{eqnarray}
\; & \; & k^{2x} \left | 
{\cos \left[y \log \left(t+{(2k-1)\over (N+1)}\right)\right] \over
\left(t+{(2k-1)\over (N+1)}\right)^{x}}
-{\cos \left[y \log \left(t+{2k \over (N+1)}\right)\right]
\over \left(t+{2k \over (N+1)}\right)^{x}} \right |
\nonumber \\
&=& 
k^{x} \left | \cos \left[y \log \left(t+{(2k-1)\over (N+1)}\right) \right]
{k^{x}\over \left(t+{(2k-1)\over (N+1)}\right)^{x}} \right .
\nonumber \\
&-& \left . \cos \left[y \log \left(t+{2k \over (N+1)}\right)\right]
{k^{x} \over \left(t+{2k \over (N+1)}\right)^{x}} \right |,
\nonumber \\
&=& k^{x} (N+1)^{x} \left | \cos \left[y \log
\left(t+{(2k-1)\over (N+1)}\right)\right]
\left({{k \over (N+1)}\over {t +{(2k-1)\over (N+1)}}}\right)^{x} \right .
\nonumber \\
&-& \left . \cos \left[y \log \left(t+{2k \over (N+1)}\right) \right]
\left({{k \over (N+1)}\over {t+{2k \over (N+1)}}}\right)^{x} \right | 
\nonumber \\
& \leq & k^{x}(N+1)^{x} \left | \cos \left[y \log \left(t+{(2k-1)\over (N+1)}\right)\right]
-\cos \left[y \log \left(t+{2k \over (N+1)}\right)\right] \right |
\nonumber \\
&=& 2 k^{x}(N+1)^{x} \left | \sin \left[{y \over 2} 
\log \left( \left(t+{(2k-1)\over (N+1)}\right)
\left(t+{2k \over (N+1)}\right)\right)\right] \right |
\nonumber \\
& \times & \left | \sin \left[{y \over 2} \log
{\left(t+{(2k-1)\over (N+1)}\right)
\over
\left(t+{2k \over (N+1)}\right)}
\right] \right |
\nonumber \\
& \leq & 
2 k^{x} (N+1)^{x} 
\left | \sin \left[{y \over 2} \log
{\left(t+{(2k-1)\over (N+1)}\right)
\over
\left(t+{2k \over (N+1)}\right)}
\right] \right | 
\nonumber \\
&=& 2 k^{x} (N+1)^{x} \left |
\sin \left[{y \over 2} \log
{(t(N+1)+2k-1)\over (t(N+1)+2k)}\right] \right |
\nonumber \\
&=& 2 k^{x} (N+1)^{x} \left | \sin \left[{y \over 2} \log
\left(1-{1 \over (t(N+1)+2k)}\right) \right] \right | 
\nonumber \\
& \le & {y k^{x}(N+1)^{x}\over (t(N+1)+2k)} ,
\label{(6.13)}
\end{eqnarray}
and this approaches $0$ as $k \rightarrow \infty$ in every compact set
\begin{equation}
{\mathbb S}_{\varepsilon,T}^{+} = \left \{s=(x+iy) \in {\mathbb S}^{+}:
\; x \in \left[{1 \over 2}+\varepsilon,1-\varepsilon \right], \;
y \in [0,T] \right \}.
\label{(6.14)}
\end{equation}
\vskip 0.3cm
\noindent
{\bf Second step}. It is now possible to prove that the $S_{L,N}(t)$ defined in (6.5) are, 
with respect to $L$, Cauchy-Ces\`aro sequences of functions according to our definition. 
Indeed, for all $t \in [1,A]$ one finds
\begin{eqnarray}
\; & \; & \Bigr | S_{L+1,N}(t)-S_{L,N}(t) \Bigr |_{(C,1)} \equiv 
\left | {1 \over M} \sum_{L=1}^{M}(S_{L+1,N}(t)-S_{L,N}(t)) \right |
\nonumber \\
&=& \left | {1 \over M} \sum_{L=1}^{M} \left \{
\sum_{k=1}^{L+1}(-1)^{k+1} \left[{1 \over \zeta(2x)}
\sum_{n=1}^{N}{1 \over n^{2x}}
{\cos \left[y \log \left(t+{k \over n}\right)\right]
\over \left(t+{k \over n}\right)^{x}}\right]
\right . \right .
\nonumber \\
&-&  \left . \left . \sum_{k=1}^{L}(-1)^{k+1}
\left[{1 \over \zeta(2x)}\sum_{n=1}^{N}
{1 \over n^{2x}}
{\cos \left[y \log \left(t+{k \over n}\right)\right]
\over \left(t+{k \over n}\right)^{x}}\right]\right \} \right |
\nonumber \\
&=& \left | {1 \over M} \sum_{L=1}^{M}
(-1)^{L} {1 \over \zeta(2x)} \sum_{n=1}^{N}
{1 \over n^{2x}}
{\cos \left[y \log \left(t+{(L+1)\over n}\right)\right] \over
\left(t+{(L+1)\over n}\right)^{x}} \right | 
\nonumber \\
&=& \left | {1 \over M} \sum_{L=1}^{[{M \over 2}]}
{1 \over \zeta(2x)} \sum_{n=1}^{N}{1 \over n^{2x}}
\left[{\cos \left[y \log \left(t+{2L \over n}\right)\right]
\over \left(t+{2L \over n}\right)^{x}} \right . \right .
\nonumber \\
&-& \left . \left . {\cos \left[y \log \left(t+{(2L+1) \over n}\right)\right]
\over \left(t+{(2L+1) \over n}\right)^{x}} \right] \right |
\nonumber \\
& \leq & {1 \over M} \sum_{L=1}^{[{M \over 2}]}
{1 \over \zeta(2x)} \sum_{n=1}^{N}{1 \over n^{2x}}
\left | {\cos \left[y \log \left(t+{2L \over n}\right)\right]
\over \left(t+{2L \over n}\right)^{x}} \right .
\nonumber \\ 
&-& \left . {\cos \left[y \log \left(t+{(2L+1) \over n}\right)\right]
\over \left(t+{(2L+1) \over n}\right)^{x}} \right |.
\label{(6.15)}
\end{eqnarray}
Now we exploit the Taylor expansion formula with remainder in the Lagrange form to find
($t_{0}^{*}$ being a point in the open interval $]t_{0},t[$)
\begin{eqnarray}
\; & \; & {\cos (y \log t)\over t^{x}}
={\cos (y \log t_{0}) \over t_{0}^{x}}
+{{\rm d}\over {\rm d}t}
\left . 
{\cos (y \log t)\over t^{x}} \right |_{t_{0}^{*}} (t-t_{0}) 
\nonumber \\
&=& {\cos (y \log t_{0}) \over t_{0}^{x}}
+{-y \sin (y \log t_{0}^{*})-x \cos (y \log t_{0}^{*})\over
(t_{0}^{*})^{x+1}}(t-t_{0}),
\label{(6.16)}
\end{eqnarray}
and hence we can write, for all $t \geq 1$ and $(x+iy) \in {\mathbb S}^{+}$,
$\alpha_{n}$ lying in the open interval $]0,1[$,
\begin{eqnarray}
\; & \; &
\left | {\cos \left[y \log \left(t+{2L \over n}\right)\right]
\over \left(t+{2L \over n}\right)^{x}} 
-{\cos \left[y \log \left(t+{(2L+1) \over n}\right)\right]
\over \left(t+{(2L+1) \over n}\right)^{x}} \right |
\nonumber \\
&=& {1 \over n} \left | 
{y \sin \left(y \log \left(t+{(2L+1+\alpha_{n})\over n}\right)\right)
\over \left(t+{(2L+1+\alpha_{n} \over n}\right)^{x+1}} \right .
\nonumber \\ 
&+& \left . {x \cos \left(y \log \left(t+{(2L+1+\alpha_{n})\over n}\right)\right) 
\over \left(t+{(2L+1+\alpha_{n})\over n}\right)^{x+1}} \;  \right |
\nonumber \\
&=& n^{x} \left |  
{y \sin \left(y \log \left(t+{(2L+1+\alpha_{n})\over n}\right)\right) +x \cos
\left(y \log \left(t+{(2L+1+\alpha_{n})\over n}\right)\right) 
\over (nt+2L+1+\alpha_{n}))^{x+1}} \right |
\nonumber \\
& \leq & (y+x) {n^{x}\over (nt+2L+1+\alpha_{n})^{x+1}}
\nonumber \\
& < & (y+x){(nt+2L+1+\alpha_{n})^{x}\over 
(nt+2L+1+\alpha_{n})^{x+1}}
\nonumber \\
& < & {(y+x) \over (nt+2L+1)},
\label{(6.17)}
\end{eqnarray}
from which we obtain the majorization
\begin{eqnarray}
\; & \; &
{1 \over \zeta(2x)} \sum_{n=1}^{N}{1 \over n^{2x}}
\left | {\cos \left[y \log \left(t+{2L \over n}\right)\right]
\over \left(t+{2L \over n}\right)^{x}} 
-{\cos \left[y \log \left(t+{(2L+1) \over n}\right)\right]
\over \left(t+{(2L+1) \over n}\right)^{x}} \right |
\nonumber \\
& \leq & {1 \over \zeta(2x)} \sum_{n=1}^{N}
{1 \over n^{2x}}{(x+y) \over (nt+2L+1)}
< {(x+y)\over 2L} {1 \over \zeta(2x)} \sum_{n=1}^{N}{1 \over n^{2x}}
\nonumber \\
& < & {(x+y)\over 2L}.
\label{(6.18)}
\end{eqnarray}
The Ces\`aro average (6.15) can be therefore majorized according to
\begin{eqnarray}
\; & \; & \Bigr | S_{L+1,N}(t)-S_{L,N}(t) \Bigr |_{(C,1)}
\nonumber \\
& \leq & {1 \over M} \sum_{L=1}^{[{M \over 2}]}
{1 \over \zeta(2x)} \sum_{n=1}^{N}{1 \over n^{2x}}
\left | {\cos \left[y \log \left(t+{2L \over n}\right)\right]
\over \left(t+{2L \over n}\right)^{x}} 
-{\cos \left[y \log \left(t+{(2L+1) \over n}\right)\right]
\over \left(t+{(2L+1) \over n}\right)^{x}} \right |
\nonumber \\
& < & {1 \over M} \sum_{L=1}^{[{M \over 2}]} {1 \over \zeta(2x)}
\sum_{n=1}^{N}{1 \over n^{2x}}
{(y+x)\over (nt+2L+1)}
\nonumber \\
& < & {(y+x)\over M} \sum_{L=1}^{[{M \over 2}]}{1 \over 2L}
\nonumber \\
& < & {(y+x)\over 2M}\log \left[{M \over 2}\right].
\label{(6.19)}
\end{eqnarray}
Therefore, for all $\varepsilon >0$ one can find a $\nu_{\varepsilon}$ such that,
for any $M > \nu_{\varepsilon}$, one achieves convergence according to our
definition (see (6.5)). It now remains to be proved that
$\left | S_{L+1,N}-S_{L,N+1} \right |_{(C,1)} < \varepsilon$.
\vskip 0.3cm
\noindent
{\bf Third step}. We point out preliminarily that
\begin{eqnarray}
\; & \; & S_{L+1,N}-S_{L,N+1} 
\nonumber \\
&=& \sum_{k=1}^{L+1}(-1)^{k+1} \left \{
{1 \over \zeta(2x)} \sum_{n=1}^{N} {1 \over n^{2x}}
{\cos \left[y \log \left(t+{k \over n}\right)\right] \over
\left(t+{k \over n}\right)^{x}} \right \}
\nonumber \\
&-& \sum_{k=1}^{L}(-1)^{k+1} \left \{ {1 \over \zeta(2x)}
\sum_{n=1}^{N+1}{1 \over n^{2x}}
{\cos \left[y \log \left(t+{k \over n}\right)\right]
\over \left(t+{k \over n}\right)^{x}} \right \}
\nonumber \\
&=& (-1)^{L} \left \{ {1 \over \zeta(2x)} \sum_{n=1}^{N}
{1 \over n^{2x}}
{\cos \left[y \log \left(t+{(L+1)\over n}\right)\right] \over
\left(t+{(L+1)\over n}\right)^{x}} \right \}
\nonumber \\
&-& \sum_{k=1}^{L}(-1)^{k+1} \left \{
{1 \over \zeta(2x)}{1 \over (N+1)^{2x}}
{\cos \left[y \log \left(t+{k \over (N+1)}\right)\right] \over
\left(t+{k \over (N+1)}\right)^{x}} \right \}
\nonumber \\
&=& {1 \over \zeta(2x)} \left \{ (-1)^{L} \sum_{n=1}^{N}
{1 \over n^{2x}} 
{\cos \left[y \log \left(t+{(L+1)\over n}\right)\right] \over
\left(t+{(L+1)\over n}\right)^{x}} \right .
\nonumber \\
&-& \left . {1 \over (N+1)^{2x}} \sum_{k=1}^{L} (-1)^{k+1}
{\cos \left[y \log \left(t+{k \over (N+1)}\right)\right] \over
\left(t+{k \over (N+1)}\right)^{x}} \right \}.
\label{(6.20)}
\end{eqnarray}
Therefore, on taking the Ces\`aro average of (6.20) we find
\begin{eqnarray}
\; & \; & \left | {1 \over M} \sum_{L=1}^{M} (S_{L+1,N}-S_{L,N+1}) \right | 
\nonumber \\
&=& \left | {1 \over M} \sum_{L=1}^{M} (-1)^{L+1}
{1 \over \zeta(2x)} \sum_{n=1}^{N} {1 \over n^{2x}}
{\cos \left[y \log \left(t+{(L+1)\over n}\right) \right] \over
\left(t+{(L+1)\over n}\right)^{x}}  \right .
\nonumber \\
&-& \left . {1 \over \zeta(2x)} {1 \over (N+1)^{2x}}
\sum_{k=1}^{L}(-1)^{k+1}
{\cos \left[y \log \left(t+{k \over (N+1)}\right)\right] \over
\left(t+{k \over (N+1)}\right)^{x}} \right |
\nonumber \\
& \leq & \left | {1 \over M} \sum_{L=1}^{[{M \over 2}]}
{1 \over \zeta(2x)}\sum_{n=1}^{N} {1 \over n^{2x}}
\left | {\cos \left[y \log \left(t+{2L \over n}\right)\right]
\over \left(t+{2L \over n}\right)^{x}} 
-{\cos \left[y \log \left(t+{(2L+1) \over n}\right)\right]
\over \left(t+{(2L+1) \over n}\right)^{x}} \right | \right .
\nonumber \\
&-& \left . {1 \over \zeta(2x)} \sum_{k=1}^{L} (-1)^{k+1}
b_{k}^{(N)}(t) \right | 
\nonumber \\
& \leq & {1 \over M} \sum_{L=1}^{[{M \over 2}]} {1 \over \zeta(2x)}
\sum_{n=1}^{N} {1 \over n^{2x}}
\left | {\cos \left[y \log \left(t+{2L \over n}\right)\right]
\over \left(t+{2L \over n}\right)^{x}} 
-{\cos \left[y \log \left(t+{(2L+1) \over n}\right)\right]
\over \left(t+{(2L+1) \over n}\right)^{x}} \right |
\nonumber \\
&+& \left | {1 \over \zeta(2x)} \sum_{k=1}^{L} (-1)^{k+1}
b_{k}^{(N)}(t) \right | 
\nonumber \\
& \leq & {(y+x)\over 2M} \log \left[{M \over 2} \right]
+{1 \over \zeta(2x)} \left | \sum_{k=1}^{L} (-1)^{k+1}
b_{k}^{(N)}(t) \right | < 2 \epsilon ,
\label{(6.21)}
\end{eqnarray}
for all $N,M > {\rm max} \left \{\mu_{\varepsilon},\nu_{\varepsilon} \right \}$. 
The majorizations (6.9), (6.19) and (6.21) hold for all $t \in [1,A]$. The resulting
convergence is therefore uniform; moreover, since the $S_{L,N}(t)$ are sequences of
equicontinuous and uniformly bounded functions, they converge to a unique 
uniformly bounded and continuous function $S$. 
This is the desired Pringsheim theorem
within the Ces\`aro framework, and its validity ensures that, in (6.4), we can write
$$
\sum_{k=1}^{\infty}(-1)^{k+1}
{\cos \left[y \log \left(1+{k \over n}\right)\right] \over
\left(1 + {k \over n} \right)^{x}}
=\lim_{t \to 1} \sum_{k=1}^{\infty}(-1)^{k+1}
f \left(t +{k \over n}\right).
$$
 
\subsection{Completion of the proof}

In light of (6.3)-(6.21), we can write that
\begin{eqnarray}
\; & \; &
\sum_{k=1}^{\infty} (-1)^{k+1} 
{\cos \left[y \log \left(1+{k \over n}\right)\right] \over 
\left(1 + {k \over n}\right)^{x}}
= \mathop{\lim} \limits_{t \to 1}
\sum_{k=1}^{\infty} (-1)^{k+1} f \left(t+{k \over n}\right)
\nonumber \\
&=& \mathop{\lim}\limits_{t\to 1}\sum\limits_{k = 1}^\infty  
{\left( { - 1} \right)^{k + 1} \textbf{E}_{1/n}^k f\left( t \right)}   
= \mathop{\lim}\limits_{t\to 1}\left\{ {1 - \sum\limits_{k = 0}^\infty  {\left( { - 1} \right)^k 
\textbf{E}_{1/n}^k } } \right\}f\left( t \right) 
\nonumber \\
&=& \mathop{\lim}\limits_{t\to 1}\left\{ 
{1 - \frac{1}{{(1 + \textbf{E}_{1/n}) }}} \right\}f\left( t \right)
\label{(6.22)}
\end{eqnarray}
From the definitions (6.2) we have that $\frac{1}{(1+\textbf{E}_h)}=\frac12 
\textbf{M}^{-1}_h$, and we can write
\begin{eqnarray}
\; & \; &
\sum\limits_{k = 1}^\infty  {\left( { - 1} \right)^{k + 1} \frac{{\cos \left[ {y\log \left( 
{1 + \frac{k}{n}} \right)} \right]}}{{\left( {1 + \frac{k}{n}} \right)^x }}}   
\nonumber \\
&=& \mathop{\lim}\limits_{t\to 1}\left\{ {1 - \frac{1}{{(1 + \textbf{E}_{1/n}) }}} \right\}
\frac{{\cos \left[ {y\log \left( t \right)} \right]}}{{t^x }} = \mathop{\lim}\limits_{t\to 1}
\left\{ {1 - \frac{1}{2}\textbf{M}_{1/n}^{ - 1} } \right\}\frac{{\cos \left[ {y\log \left( t \right)} 
\right]}}{{t^x }}
\nonumber \\
& \equiv & \mathop{\lim}\limits_{t\to 1}\varphi_n\left(t\right). 
\label{(6.23)}
\end{eqnarray}
By applying to the last equality the operator $\textbf{M}_{1/n}$, 
whose inverse is denoted by $\textbf{M}_{1/n}^{-1}$, we find
\begin{eqnarray}
\; & \; &
\mathop{\lim}\limits_{t\to 1}\textbf{M}_{1/n} \varphi_n \left( t \right) = \mathop{\lim}
\limits_{t\to 1}\textbf{M}_{1/n} \left\{ {1 - \frac{1}{2}\textbf{M}_{1/n}^{ - 1} } \right\}
\frac{{\cos \left[ {y\log \left( {t} \right)} \right]}}{{t^x }}  
\nonumber \\
&=& \mathop{\lim}\limits_{t\to 1}\left\{\textbf{M}_{1/n} \frac{{\cos \left[ {y\log \left( {t} \right)} 
\right]}}{{t^x }} - \frac{1}{2}\frac{{\cos \left[ {y\log \left( {t} \right)} \right]}}
{{t^x }}\right\}  
\nonumber \\
&=& \frac{1}{2}\mathop{\lim}\limits_{t\to 1}\left\{ {\frac{{\cos \left[ {y\log \left( {t} \right)} 
\right]}}{{t^x }} + \frac{{\cos \left[ {y\log \left( {t + \frac{1}{n}} \right)} \right]}}
{{\left( {t + \frac{1}{n}} \right)^x }}} \right\} - \frac{1}{2}\mathop{\lim}\limits_{t\to 1}
\frac{{\cos \left[ {y\log \left( {t} \right)} \right]}}{{t^x }}  
\nonumber \\
&=& \frac{1}{2}\mathop{\lim}\limits_{t\to 1}\frac{{\cos \left[ {y\log \left( {t + \frac{1}{n}} 
\right)} \right]}}{{\left( {t + \frac{1}{n}} \right)^x }} \leq \frac{1}{2}\mathop{\lim}\limits_{t\to 1}
\left| {\frac{{\cos \left[ {y\log \left( {t + \frac{1}{n}} \right)} \right]}}
{{\left( {t + \frac{1}{n}} \right)^x }}} \right| < \frac{1}{2},  
\label{(6.24)}
\end{eqnarray}
for any $\left(x+iy\right)\in\mathbb S^{+}$, obtaining therefore 
the majorization $\mathop{\lim}\limits_{t\to 1}
\textbf{M}_{1/n}\varphi_n\left(t\right)<\frac12$. Upon remarking  that, for any constant $c\in\mathbb R$ 
we have $\textbf{M}^{-1}_{1/n} c=c+\chi_{n} \left(t\right)$, where the $\chi_{n}$ 
are real-valued functions with vanishing mean value, 
and by applying to the first and the last member of the previous equation the 
operator $\textbf{M}^{-1}_{1/n}$, we can write 
\begin{equation}
{\mathop{\lim}\limits_{t\to 1}\varphi_n \left( t \right) = \mathop{\lim}\limits_{t\to 1}
\sum\limits_{k = 1}^\infty  {\left( { - 1} \right)^{k + 1} \frac{{\cos \left[ {y\log \left( 
{t + \frac{k}{n}} \right)} \right]}}{{\left( {t + \frac{k}{n}} \right)^x }}}  
< \frac{1}{2} + \mathop{\lim}\limits_{t \to 1} \chi_{n}(t)}, 
\label{(6.25)}
\end{equation}
where $\chi_{n}$ is such that
$ 
\textbf{M}_{1/n}\chi_{n}(t)=0.  
$
Since $\chi_{n}$ has vanishing mean value for any $n\in\mathbb N$, we have necessarily that 
\begin{equation}
{\textbf{M}_{1/n} \chi_{n} \left( t \right) = \frac{1}{2}\left\{ {\chi_{n}(t) 
+ \chi_{n} \left({t + \frac{1}{n}} \right)} \right\} = 0} 
\Longrightarrow {\chi_{n}(t) 
=  - \chi_{n} \left({t + \frac{1}{n}} \right)}.
\label{(6.26)}
\end{equation}
This result must hold for any $t \in \mathbb R$ and $n \in \mathbb N$. 

At this stage, it is clear that we need to know the limit as $t$ approaches $1$ of
the solutions of Eq. (6.26). Our findings are presented in Lemmas $3$ and $4$ below. 
\vskip 0.3cm
\noindent
{\bf Lemma 3}. Any infinite-dimensional vector space of functions 
dense in the set of solutions of the functional equation (6.26)
has a basis $\left\{\psi_n^m\right\}_{m\in\mathbb N}$ consisting of anti-periodic functions as well.
\vskip 0.3cm
\noindent
{\it Proof}. By virtue of the hypothesis of density, for any $\varepsilon>0$, we can find 
a sequence $\left\{b_m\right\}_{m\in\mathbb N}$ such that
\begin{equation}
{\left| {\sum\limits_m {b_m \psi _n^m \left( t \right)}  - \chi _n \left( t \right)} \right| 
< \varepsilon}, \; \;   
{\left| {\sum\limits_m {b_m \psi _n^m \left( 
{t + \frac{1}{n}} \right)}  - \chi _n \left( {t + \frac{1}{n}} \right)} \right| < \varepsilon }, 
\label{(6.27)}
\end{equation}
for any $n\in\mathbb N$. From the hypothesis $\chi _n \left( {t + \frac{1}{n}} \right) 
=  - \chi _n \left( t \right)$ we can write
\begin{eqnarray}
\; & \; &
\biggr | \sum\limits_m {b_m \left\{ { \psi _n^m \left( t \right) 
+  \psi _n^m \left( {t + \frac{1}{n}} \right)} \right\}}  \biggr |  
\le \biggr | \sum\limits_m {b_m  \psi _n^m \left( t \right)}  - \chi _n \left( t \right) \biggr |
\nonumber \\ 
&+& \biggr | \chi_{n} \left(t \right) + \chi_{n} \left({t + \frac{1}{n}}\right) \biggr|
\nonumber \\ 
&+& \biggr| \sum\limits_m {b_m  \psi _n^m \left( {t + \frac{1}{n}} \right)}  
- \chi _n \left( {t + \frac{1}{n}} \right) \biggr| < 2\varepsilon   
\label{(6.28)}
\end{eqnarray}
for any choice of sequence $\left\{b_m\right\}_{m\in\mathbb N}$, hence we have that
\begin{equation}
{\sum\limits_m {b_m \left\{ { \psi _n^m \left( t \right) 
+  \psi _n^m \left( {t + \frac{1}{n}} \right)} \right\}}  
= 0 \; {\rm for \; any \; choice \; of}} \; {\left\{ {b_m } \right\}_{m \in\mathbb N} },  
\label{(6.29)}
\end{equation}
from which it follows that the vectors $\psi _n^m \left( t \right)$ and 
$\psi _n^m \left( {t + \frac{1}{n}} \right)$ depend on each other,
and in particular $ \psi _n^m \left( t \right) 
+  \psi _n^m \left( {t + \frac{1}{n}} \right)=0$, finding therefore
\begin{equation}
{\psi _n^m \left( t \right) =  - \psi _n^m \left( {t + \frac{1}{n}} \right)} 
\; {\rm for \; any} \; {m \in\mathbb N},  
\label{(6.30)}
\end{equation}
i.e. is the thesis.

The next step of our analysis is as follows. 
\vskip 0.3cm
\noindent
{\bf Lemma 4 (All functions $\chi_{n}$ in Eq. (6.25) have vanishing limit as
$t$ approaches $1$)}. For any $n\in\mathbb N$ we have the limit
\begin{equation}
\mathop {\lim }\limits_{t \to 1} \chi_n \left( {t + \frac{1}{n}} \right) 
=  - \mathop {\lim }\limits_{t \to 1} \chi_n \left( t \right) = 0.
\label{(6.31)}
\end{equation}
\vskip 0.3cm
\noindent
{\it Proof}. On taking into account the \textsc{Gram-Schmidt Theorem}, 
and exploiting the basis 
$$
\left\{\psi_n^m\right\}_{m\in\mathbb N}
$$ 
we can find one and only one 
ortho-normalized basis $\left\{\tilde\varphi_n^m\right\}_{m\in\mathbb N}$ such that
(on denoting by $\delta^{kl}$ the Kronecker symbol, equal to $1$ if 
$k=l$ and equal to $0$ if $k \not =l$)
\begin{equation}
{\left( {\tilde \varphi_{n}^{k} \left( t \right),\tilde \varphi_{n}^{l} \left( t \right)} \right) 
= \delta^{kl} } \; {\rm for} \; {k,l \in\mathbb N},  
\label{(6.32)}
\end{equation}
for any $t\in {\mathbb R}^{+}$. Now, by choosing a positive real number $A$ in the
open interval $]1,\infty[$ (see below), the ortho-normal basis can be built
for $t \in [1,A]$, finding therefore 
\begin{equation}
{\tilde\varphi_{n}^{m} \left( t \right) = a_{nm} \sin \left( {n m \pi t} \right)}   
\label{(6.33)}
\end{equation}
for any $n,m \in {\mathbb N}$, where we have set
\begin{equation}
a_{nm} = {1 \over \sqrt{\int_{1}^{A}\sin^{2} (mn \pi t)dt}},
\label{(6.34)}
\end{equation}
whose denominator never vanishes if $A > 1$. Hence we can write
\begin{equation}
{\underbrace {\mathop {\lim }\limits_{t \to 1} \sum\limits_{m = 1}^\infty  {a_{nm}b_m \sin \left[ 
{nm\pi \left( {t + \frac{1}{n}} \right)} \right]} }_{\mathop {\lim }\limits_{t \to 1} \chi _n \left( 
{t + \frac{1}{n}} \right)} = - \underbrace { \mathop {\lim }\limits_{t \to 1} \sum\limits_{m = 1}^\infty  
{a_{nm} b_m \sin \left( {nm\pi t} \right)} }_{ \mathop {\lim }\limits_{t \to 1} \chi _n \left( t \right)}}   
\label{(6.35)}
\end{equation}
for any $n \in {\mathbb N}$, obtaining therefore the thesis.\footnote{We note incidentally that, 
from Eq. (6.26), it follows that
$$
\chi_{n}(t)=(-1)^{k}\chi_{n} \left(t+{k \over n}\right) \; \; 
\forall k=1,2,...,
$$
which is of course satisfied by the functions in Eq. (6.33).}

Having proved that $\chi_{n}$ is vanishing for all $n$ as $t$ approaches $1$,
taking into account the \textsc{Pringsheim Theorem} about order of summation 
exchange proved in Appendix B, we can write (see (3.9) and (6.4)) 
\begin{eqnarray}
\; & \; &
\sum\limits_{k = 1}^\infty  {\left( { - 1} \right)^{k + 1} Z_{k + 1} \left( {x,y} \right)}
\nonumber \\  
&=& \mathop{\lim}\limits_{t\to 1}\sum\limits_{k = 1}^\infty  {\left( { - 1} \right)^{k + 1} } 
\left\{ {\frac{1}{{\zeta \left( {2x} \right)}}\sum\limits_{n = 1}^\infty  {\frac{1}{{n^{2x} }}} 
\frac{{\cos \left[ {y\log \left( {t + \frac{k}{n}} \right)} \right]}}
{{\left( {t + \frac{k}{n}} \right)^x }}} \right\}  
\nonumber \\
&=& \mathop{\lim}\limits_{t\to 1}\frac{1}{{\zeta \left( {2x} \right)}}\sum\limits_{n = 1}^\infty  
{\frac{1}{{n^{2x} }}} \left\{ {\sum\limits_{k = 1}^\infty  {\left( { - 1} \right)^{k + 1} } 
\frac{{\cos \left[ {y\log \left( {t + \frac{k}{n}} \right)} \right]}}
{{\left( {t + \frac{k}{n}} \right)^x }}} \right\} 
\nonumber \\
&=& \mathop{\lim}\limits_{t\to 1}\frac{1}
{{\zeta \left( {2x} \right)}}\sum\limits_{n = 1}^\infty  {\frac{1}{{n^{2x} }}
\varphi_n \left( t \right)}  < \frac{1}{2} , 
\label{(6.36)}
\end{eqnarray}
obtaining therefore the desired contradiction. 

\section{Concluding remarks}

At the risk of slight repetitions, we find it appropriate to summarize the original parts of 
our classical analysis approach to the Riemann hypothesis as follows.
\vskip 0.3cm
\noindent
(i) A study of the double series (3.1) in ${\mathbb S}^{+}$, with the associated
zeros' functions defined in (3.9).
\vskip 0.3cm
\noindent
(ii) Derivation of the fundamental identity (4.1), with its maximal
extension (4.3).
\vskip 0.3cm
\noindent
(iii) Proof in Appendix B of the Pringsheim convergence (cf. 
Refs. \cite{rif4,Brom,Hardy,Lima}), which is necessary to make sure that the steps in
(i) and (ii) are meaningful.
\vskip 0.3cm
\noindent
(iv) Characterization (5.1) of non-trivial zeros of the Riemann $\zeta$-function, and
invalidation (6.36) of this condition with the help of finite-difference operators defined
in (6.2) and of the detailed sub-section $6.1$.
\vskip 0.3cm
\noindent
As far as we can see, what remains to be done within our classical 
analysis framework is to prove an uniqueness theorem for the solution
of the functional equation (6.26). Our aim was to develop a general classical
analysis framework where the approach to proving the Riemann hypothesis relies
upon the explicit proof of uniform convergence of double series. Our equivalent
of the Riemann hypothesis is, as far as we can see, falsifiable with the help 
of classical analysis only.

Nevertheless, we are aware of the power of the methods of modern analysis, that
are able to relate all equivalent formulations found so far \cite{Broughan}.
The desired proof would provide
therefore fundamental contributions, at the same time, to all branches
of mathematics where the Riemann hypothesis plays a role
(e.g. \cite{KA,Bombieri,Langl,Bristol}).
 
\begin{appendix}

\section{Convergence of monotonic alternating double series and Pringsheim's theorem}
\setcounter{equation}{0}

Before analyzing the properties of Riemann's $\zeta$-function in ${\mathbb S}^+$, we need to 
get rid of any embarrassment about which Criterion of Convergence is used to give meaning to 
all series we have met in our investigation. In particular, we will need 
to study the conditions under which the sums
used in section 3 and defined in ${\mathbb S}^+$ are meaningful. 

Given an arbitrary function
\begin{equation}
a:{\mathbb N} \times {\mathbb N} \longrightarrow a_{ij} \in {\mathbb C}
\label{(A.1)}
\end{equation}
its image is the infinite {\it lattice} 
\begin{equation}
{\cal L} =   \left\{\begin{array}{*{20}c}
{a_{11} } & {a_{12} } &  \cdots  & {a_{1j} } &  \cdots   \\
{a_{21} } & {a_{22} } &  \cdots  & {a_{2j} } &  \cdots   \\
\vdots  &  \vdots  &  \ddots  &  \vdots  & {}  \\
{a_{i1} } & {a_{i2} } &  \cdots  & {a_{ij} } &  \cdots   \\
\vdots  &  \vdots  & {} &  \vdots  &  \ddots   \\
\end{array}\right\} .
\label{(A.2)}
\end{equation}
We can ask whether the sum of {\it all} its terms, that we call {\it double series}, 
does exist and, in the affirmative case, whether it is finite or infinite. 
In order to specify what definition of convergence we adopt, we have to define 
different subsets in the lattice ${\cal L}$. The Pringsheim Region \cite{rif4} 
\begin{equation}
{\cal R}_{pq} = \left\{a_{ij}\right\}_{i\le p,j\le q}
\label{(A.3)}
\end{equation}
is a connected rectangular subset of the lattice containing any element 
$a_{ij}$ with $i\le p$ and $j\le q$:  
\begin{equation}
{\cal R}_{pq} = \left\{ {\begin{array}{*{20}c}
{a_{11} } & {a_{12} } &  \cdots  & {a_{1q} }   \\
{a_{21} } & {a_{22} } &  \cdots  & {a_{2q} }    \\
\vdots  &  \vdots  &  \ddots  &  \vdots    \\
{a_{p1} } & {a_{p2} } &  \cdots  & {a_{pq} }   \\
\end{array}} \right\}.
\label{(A.4)}
\end{equation}
On it we define a {\it partial sum} on a Pringsheim set \cite{rif4} as
\begin{equation}
\begin{array}{*{20}c}
{S_{pq} = \sum_{i \le p \hfill \atop 
 j \le q \hfill} {a_{ij} } } & {\rm{or}} & {S_{pq}  =
\sum_{\left\{ {{\cal R}_{pq} } \right\}} {a_{ij} } } . \\
\end{array}
\label{(A.5)}
\end{equation}
Thus, we can introduce the {\it Pringsheim Convergence} criterion \cite{rif4} as
\begin{equation}
\begin{array}{*{20}c}
{\forall \varepsilon  > 0} & {\exists \nu _\varepsilon   \in {\mathbb N}} & {\exists S   
\in {\mathbb C}} & : & {\forall p,q > \nu _\varepsilon} &  \Rightarrow  
& {\left| {S_{pq}  - S} \right| < \varepsilon },  
\end{array}
\label{(A.6)}
\end{equation}
where $i$ and $j$ run from $1$ through $p$ and $q$ independently\footnote{The 
statement that the indices $i,j$ run \textsl{at the 
same time but independently} means that none of them is constrained by a particular algorithm, 
but both are free to run in ${\mathbb N}$ independently one of the other. 
The partial sum of a double series in the Pringsheim convergence 
splits therefore ahead as follows:
\begin{equation}
\begin{array}{*{20}c}
{\sum_{i=1}^N {\sum_{j=1}^{M} a_{ij} } } 
& \mathop{\longrightarrow}\limits^{\rm{split \; ahead}} 
& {\sum_{i=1}^{N + X_i } {\sum_{j=1}^{M + X_j } {a_{ij} } } }.  \\
\end{array}
\label{(A.7)}
\end{equation}
}. We refer to it with the concise notation 
$\sum_{\left\{ {\cal R} \right\}} {a_{ij} }  < \infty $, 
where with ${\cal R}$ we consider the whole class of rectangles ${\cal R}_{pq}$. 
Two other very important Pringsheim Sets are the {\it columns} ${\cal R}_{p\bar h}$ and the 
{\it rows} ${\cal R}_{\bar h q}$, where $\bar h$ is an index taking a fixed value. 
We can now define the {\it row partial sum} and the {\it column partial sum} as the 
sums\footnote{In the next formulas of partial sums the fixed indices are $i$ or $j$ and will not 
be over-scored by a bar.}
\begin{equation}
\begin{array}{*{20}c}
{s_{ i,q + l}  = \sum \limits_{j = 1}^{q + l} {a_{ij} } },  & \; 
{s_{p+m, j} = \sum \limits_{i=1}^{p+m} {a_{ij} } },  \\
\end{array}
\label{(A.8)}
\end{equation}
and then introduce the concepts of {\it column-convergence} 
\begin{equation}
\begin{array}{*{20}c}
{\forall \varepsilon  > 0} & {\exists \nu _\varepsilon   \in {\mathbb N}} & {\exists c_{\bar h}   
\in {\mathbb C}} & : & {\forall p > \nu _\varepsilon} &  \Rightarrow  
& {\left| {s_{p\bar h}  - c_{\bar h}} \right| < \varepsilon },  
\end{array}
\label{(A.9)}
\end{equation}
and {\it row-convergence} 
\begin{equation}
\begin{array}{*{20}c}
{\forall \varepsilon  > 0} & {\exists \nu _\varepsilon   \in {\mathbb N}} & {\exists r_{\bar h}   
\in {\mathbb C}} & : & {\forall q > \nu _\varepsilon} &  \Rightarrow  
& {\left| {s_{\bar h q}  - r_{\bar h}} \right| < \varepsilon}.  
\end{array}
\label{(A.10)}
\end{equation}
After the introduction of these definitions, Pringsheim was able to prove the following \cite{rif4}
\vskip 0.3cm
\noindent
{\bf Theorem 3 (Pringsheim)}. If the double-series 
$$
\sum\limits_{i,j=1}^{\infty} a_{ij}  		
$$
is Pringsheim-, column- and row-convergent one can exchange the order of summation, i.e. 
\begin{equation}
\sum\limits_{i = 1}^\infty  {\left\{ {\sum\limits_{j = 1}^\infty  {a_{ij} } } \right\}}  
= \sum\limits_{j = 1}^\infty  {\left\{ {\sum\limits_{i = 1}^\infty  {a_{ij} } } \right\}}.
\label{(A.11)}
\end{equation}
\vskip 0.3cm
\noindent
{\it Proof}. From the hypothesis of column- and row-convergence, for 
any $\varepsilon > 0$ we can find a $\nu_{\varepsilon} \in {\mathbb N}$ such that,
for any $p,q > \nu_{\varepsilon}$, we have 
\begin{equation}
\left | \sum_{j=1}^{q}a_{ij}-r_{i} \right | = \left | s_{i,q}-r_{i} \right |
< \varepsilon, \; \;
\left | \sum_{i=1}^{p}a_{ij}-c_{j} \right | 
=\left | s_{p,j} -c_{j} \right | < \varepsilon.
\label{(A.12)}
\end{equation}
Now we can see that, by increasing in a suitable way the limit of summation we can
reduce in a convenient way the $\varepsilon$ upper bound; in particular, we can find
two numbers $l=l(p) \in {\mathbb N}$ and $m=m(q) \in {\mathbb N}$ such that
\begin{eqnarray}
\; & \; & 
\left| \sum_{j=1}^{q+l(p)} a_{ij}   - r_{i}  \right| = \left| s _{i,q+l}-r_{i}\right| 
< {\varepsilon \over \zeta(2) p^{2}} 
\leq {\varepsilon \over \zeta(2)i^{2}} \; 
\forall i \in \left \{ 1,...,p \right \} , 
\nonumber \\
& \; & \left| \sum_{i=1}^{p+m(q)} a_{ij}   - c_{j}  \right| 
= \left| s _{p+m,j}-c_{j}\right|  
\nonumber \\
& < & {\varepsilon \over \zeta(2) q^{2}}
\leq {\varepsilon \over \zeta(2) j^{2}} \;
\forall j \in \left \{ 1,...,q \right \},
\label{(A.13)}
\end{eqnarray} 
hence we have
\begin{equation}
\left| {\sum\limits_{i = 1}^{p + m} {\left( {s_{ i,q + l}  - r_{ i} } \right)} } \right| 
\leq \sum\limits_{i = 1}^{p + m} {\left| {s_{i,q + l}  - r_{ i} } \right|}  
< \sum\limits_{i = 1}^{p + m} 
{\frac{\varepsilon }{{\zeta \left( 2 \right)i^2 }}}  
< \varepsilon ,
\label{(A.14)}
\end{equation}
\begin{equation}
\left| {\sum\limits_{j = 1}^{q+l} {\left( {s_{ p+m,j}  - c_{j} } \right)} } \right| 
\leq \sum\limits_{j=1}^{q+l} {\left| {s_{p+m,j}  - c_{j} } \right|}  
< \sum\limits_{j=1}^{q+l} 
{\frac{\varepsilon }{{\zeta \left( 2 \right)j^2 }}}  
< \varepsilon .
\label{(A.15)}
\end{equation}
From the hypothesis of Pringsheim convergence we can say that, for any $\varepsilon>0$, we can find a 
$\nu_{\varepsilon}\in {\mathbb N}$ and a $S \in {\mathbb C}$ 
such that, for any $p,q>\nu_{\varepsilon}$, we have 
$\left| {S_{pq}  - S} \right| < \varepsilon $; moreover, for 
any $m$ and $l$ in ${\mathbb N}$, we can write
\begin{equation}
\begin{array}{*{20}c}
{\left| {\sum\limits_{i = 1}^{p + m} {\sum\limits_{j = 1}^{q + l} {a_{ij} } }  - S} \right| 
= \left| {S_{p+m,q+l}  - S} \right| < \varepsilon }, & {\forall p,q > \nu _\varepsilon  }, 
& {\forall l,m \in {\mathbb N}}.  
\end{array}
\label{(A.16)}
\end{equation}
Now, in light of
\begin{equation}
\begin{array}{*{20}c}
{\sum\limits_{i = 1}^{p + m} {\sum\limits_{j = 1}^{q + l} {a_{ij} } }  - S 
= \sum\limits_{i = 1}^{p + m} {\left( {s_{i,q + l}  - r_{i} } \right)}  
+ \sum\limits_{i = 1}^{p + m} {r_{i} }  - S}, & {\forall p,q \in {\mathbb N}}, 
& {\forall l,m \in {\mathbb N}},  \\
\end{array}
\label{(A.17)}
\end{equation}
bearing in mind (A.14) and (A.15), we can write
\begin{eqnarray}
\; & \; &
\left| {\sum\limits_{i = 1}^{p + m} {r_{ i} }  - S} \right| = \left| 
{\left( {\sum\limits_{i = 1}^{p + m} {\sum\limits_{j = 1}^{q + l} {a_{ij} } }  - S} \right)
-\sum\limits_{i = 1}^{p + m} {\left( {s_{i,q + l}  - r_{ i} } \right)}   } \right|     
\nonumber \\
& < & \left| {\sum\limits_{i = 1}^{p + m} {\sum\limits_{j = 1}^{q + l} {a_{ij} } }  - S} \right| 
+\left| {\sum\limits_{i = 1}^{p + m} {\left( {s_{i,q + l}  - r_{ i} } \right)} } \right|  
\nonumber \\
& < & 2\varepsilon , \; \forall p,q > \nu_{\varepsilon}, \; \forall l,m \in {\mathbb N},
\label{(A.18)}
\end{eqnarray}
hence, by adding and subtracting partial sums of columns and rows and exploiting (A.13), 
(A.14) and (A.17) we have eventually
\begin{eqnarray}
\; & \; & 
\left| {\sum\limits_{i = 1}^{p + m} {\sum\limits_{j = 1}^{q + l} {a_{ij} } }  
- \sum\limits_{j=1}^{q + l} {\sum\limits_{i = 1}^{p + m} {a_{ij} } } } \right|  
\nonumber \\
&=& \left| \sum_{i=1}^{p+m} \left(\sum_{j=1}^{q+l}a_{ij}-r_{i}\right)
-\sum_{j=1}^{q+l}\left(\sum_{i=1}^{p+m}a_{ij}-c_{j}\right) \right .
\nonumber \\
&+& \left . \left(\sum_{i=1}^{p+m}r_{i}-S \right)
-\left(\sum_{j=1}^{q+l}c_{j}-S \right) \right |
\nonumber \\ 
& < & \sum\limits_{i = 1}^{p + m} {\left| 
{\sum\limits_{j = 1}^{q + l} {a_{ij} }  - r_{ i} } \right|}  
+ \sum\limits_{j = 1}^{q + l} {\left| {\sum\limits_{i = 1}^{p + m} {a_{ij} }  - c_{ j} } \right|}  
+ \left| {\sum\limits_{i = 1}^{p + m} {r_{ i} }  - S} \right| 
\nonumber \\
&+& \left| {\sum\limits_{j = 1}^{q + l} 
{c_{ j} }  - S} \right| 
\nonumber \\
& < & 4\varepsilon .   
\label{(A.19)}
\end{eqnarray}
This result being true for arbitrary choice of $l,m\in {\mathbb N}$, 
we obtain the thesis when $l,m\to\infty$. Q.E.D.

\section{Pringsheim convergence of particular double series in ${\mathbb S}^+$}
\setcounter{equation}{0}

We have to prove that the series we have introduced in Sec. 3
\begin{equation}
\begin{array}{*{20}c}
{\sum\limits_{n_1  \ne n_2 } {\frac{{\left( { - 1} \right)^{n_1  + n_2 } }}{{n_1^s n_2^{\bar s} }}} } 
& {;} & {\sum\limits_{n = 1}^\infty  {\frac{{\cos \left[ {y\log \left( {1 + \frac{k}{n}} \right)} \right]}}
{{\left[ {n\left( {n + k} \right)} \right]^x }}} } \; {\rm and} \; {\sum\limits_{k = 1}^\infty  
{\left( { - 1} \right)^{k + 1} {\sum\limits_{n = 1}^\infty  {\frac{{\cos \left[ {y\log
\left( {1 + \frac{k}{n}} \right)} \right]}}{{\left[ {n\left( {n + k} \right)} \right]^x }}} }} }  \\
\end{array}
\label{(B.1)}
\end{equation}
are meaningful. For this purpose, we begin by remarking that, since the $Z_{k+1}$ functions are
uniformly bounded and equicontinuous for all $\varepsilon >0$, when the Ascoli-Arzel\`a theorem is
applied to every compact subset\footnote{Since ${\mathbb S}^{+}$ is a subset of the complex plane
where the Hausdorff axiom holds, for all $\varepsilon>0$ we find an open set 
${\cal O}_{\varepsilon}$ such that, with the notation of Eq. (1.10), the set
$\Bigr[{\mathbb S}_{\varepsilon,T}^{+}\Bigr]$ is properly included in 
${\cal O}_{\varepsilon}$. On the other hand, one has
$$
\Bigr[{\mathbb S}_{\varepsilon,T}\Bigr] \subset {\cal O}_{\varepsilon} 
\subset {\mathbb S}^{+}, \;
\forall \varepsilon >0, \; \forall T >0.
$$} 
of $S_{\varepsilon}^{+}$, their sequence converges uniformly to
an uniformly bounded and equicontinuous function. Thus, if the Pringsheim criterion (see below) holds
pointwise for numerical double series \cite{rif4,Brom,Lima}, it holds also for double series
of uniformly bounded and equicontinuous functions.

As we know from (3.2), the first double series in (B.1) can be re-written in the form
\begin{eqnarray}
\sum\limits_{n_1  \ne n_2 } {\frac{{\left( { - 1} \right)^{n_1  + n_2 } }}{{n_1^s n_2^{\bar s} }}}  
&=& \sum\limits_{n_1  \ne n_2 } {\frac{{\left( { - 1} \right)^{n_1  + n_2 } }}{{\left( {n_1 n_2 } 
\right)^x }}\cos \left( {y\log \frac{{n_2 }}{{n_1 }}} \right)}  
\nonumber \\
&+& i\sum\limits_{n_1  \ne n_2 } 
{\frac{{\left( { - 1} \right)^{n_1  + n_2 } }}{{\left( {n_1 n_2 } \right)^x }}\sin \left( 
{y\log \frac{{n_2 }}{{n_1 }}} \right)},
\label{(B.2)} 
\end{eqnarray}
hence it is sufficient to prove that
\begin{equation}
\begin{array}{*{20}c}
{\sum\limits_{n_1  \ne n_2 } {\frac{{\left( { - 1} \right)^{n_1  + n_2 } }}{{\left( {n_1 n_2 } 
\right)^x }}\cos \left( {y\log \frac{{n_2 }}{{n_1 }}} \right)}  } 
& {\rm and} & {\sum\limits_{n_1  \ne n_2 } 
{\frac{{\left( { - 1} \right)^{n_1  + n_2 } }}{{\left( {n_1 n_2 } \right)^x }}\sin \left( {y\log 
\frac{{n_2 }}{{n_1 }}} \right)}}  \\
\end{array}
\label{(B.3)}
\end{equation}
$\forall x+iy \in {\mathbb S}^{+}$
are both Pringsheim-convergent. We prove the Pringsheim-convergence of the former. Before going 
ahead we define the partial sums with respect to $n_1$ and $n_2$:
\begin{eqnarray}
\; & \; &
{s_{n_1 ,p + l}  =  \sum_{n_2  = 1}^{p + l} {
\frac{{\left( { - 1} \right)^{n_1  + n_2 } }}{{\left( {n_1 n_2 } \right)^x }}\cos \left( 
{y\log \frac{{n_2 }}{{n_1 }}} \right)} },
\nonumber \\ 
& \; & {s_{p+m,n_2 }  = 
\sum_{n_1  = 1}^{p + m} {\frac{{\left( { - 1} \right)^{n_1  + n_2 } }}
{{\left( {n_1 n_2 } \right)^x }}\cos \left( {y\log \frac{{n_2 }}{{n_1 }}} \right)} },  
\label{(B.4)}
\end{eqnarray}
and {\it row} and {\it column} sums as
\begin{eqnarray}
\; & \; &
{r_{n_2 }   = \frac{(-1)^{n_{2}}}{{n_2^x }}\sum_{n_1  = 1}^\infty  
{\frac{{\left( { - 1} \right)^{n_1 } }}{{n_1^x }}\cos \left( {y\log \frac{{n_2 }}{{n_1 }}} 
\right)}}, 
\nonumber \\ 
& \; & {c_{n_1 }  = \frac{(-1)^{n_{1}}}{{n_1^x }}
\sum_{n_2  = 1}^\infty  {\frac{{\left( { - 1} \right)^{n_2 } }}{{n_2^x }}
\cos \left( {y\log \frac{{n_2 }}{{n_1 }}} \right)}  }.  
\label{(B.5)}
\end{eqnarray}
Now, in order to go ahead it is necessary to prove the following
\vskip 0.3cm
\noindent
{\bf Lemma 5}. For any $x+iy \in {\mathbb S}^+$ we have that 
\cite{Smyrlis} (here $\alpha \in \{ 1,2 \}$)
\begin{eqnarray}
\; & \; & F_{\alpha}(x,y)  = \sum_{n_{\alpha} = 1}^\infty
{(-1)^{n_{\alpha}}\over n_{\alpha}^{x}}
\cos \left( {y\log n_{\alpha}} \right)  < \infty , 
\nonumber \\
& \; & G_{\alpha}(x,y) = \sum_{n_{\alpha}=1}^{\infty}
{(-1)^{n_{\alpha}}\over n_{\alpha}^{x}}
\sin \left(y \log n_{\alpha}\right) < \infty.
\label{(B.6)}
\end{eqnarray}
\vskip 0.3cm
\noindent
{\it Proof}. We can observe that, by separating even and odd parts of the series, we have
\begin{equation}
F_{\alpha} \left(x,y \right) = \mathop {\lim}\limits_{k \to \infty} \sum_{m = 1}^k {\left\{ 
{\frac{{\cos \left[ {y\log \left( {2m + 1} \right)} \right]}}{{\left( {2m + 1} \right)^x }} 
- \frac{{\cos \left[ {y\log \left( {2m + 2} \right)} \right]}}
{{\left( {2m + 2} \right)^x }}} \right\}}. 
\label{(B.7)}
\end{equation}
The difference in curly brackets in (B.7) can be written as follows:
\begin{eqnarray}
\; & \; &
\left| {\frac{{\cos \left[ {y\log \left( {2m + 1} \right)} \right]}}{{\left( {2m + 1} \right)^x }} 
- \frac{{\cos \left[ {y\log \left( {2m + 2} \right)} \right]}}{{\left( {2m + 2} \right)^x }}} \right|   
\nonumber \\
&=& \left| {\frac{{\left( {2m + 2} \right)^x \cos \left[ {y\log \left( {2m + 1} \right)} \right] 
- \left( {2m + 1} \right)^x \cos \left[ {y\log \left( {2m + 2} \right)} \right]}}{{\left( 
{2m + 1} \right)^x \left( {2m + 2} \right)^x }}} \right|   
\nonumber \\
&=& \left| {(2m+2)^x \left\{ {\cos \left[ {y\log \left( {2m + 1} \right)} 
\right] - \cos \left[ {y\log \left( {2m + 2} \right)} \right]} \right\}
\over (2m+1)^{x} (2m+2)^{x}} \right . 
\nonumber \\
&+& \left .  {\left\{ (2m+2)^{x}  - (2m+1)^{x}  \right\}
\cos \left[{y\ln \left( {2m + 2} \right)} \right] \over
(2m+1)^{x}(2m+2)^{x}} \right| .
\label{(B.8)}
\end{eqnarray}

Then, taking into account the {\it Integral Mean Value Theorem} we can find a 
$\lambda \in\left(0,1\right)$ such that\footnote{Thanks to the Integral Mean Value Theorem we can find 
a $\lambda \in\left(0,1\right)$ such that 
\begin{equation}
\begin{array}{l}
- \int\limits_{\left( {2k + 1} \right)}^{\left( {2k + 2} \right)} {\frac{{\sin \left( 
{y\log t} \right)}}{t}dt}  =  - \int\limits_{\log \left( {2k + 1} \right)}^{\log \left( {2k + 2} \right)} 
{\sin \left( {y\log t} \right)d\left( {\log t} \right)}  =  - \frac{1}{y}\int\limits_{y\log 
\left( {2k + 1} \right)}^{y\log \left( {2k + 2} \right)} {\sin \left( {y\log t} \right)
d\left( {y\log t} \right)}  \\ 
\\ 
= \frac{1}{y}\left\{ {\cos \left[ {y\log \left( {2k + 2} \right)} \right] - \cos \left[ 
{y\log \left( {2k + 1} \right)} \right]} \right\} = \frac{{\sin \left[ {y\log \left( 
{2k + 1 + \lambda } \right)} \right]}}{{\left( {2k + 1 + \lambda } \right)}} \\ 
\end{array}
\label{(B.9)}
\end{equation}
Moreover, we can apply it once more to find
\begin{equation}
\int\limits_{\left( {2k + 1} \right)}^{\left( {2k + 2} \right)} {xt^{x - 1} dt}  
= \left. {t^x } \right|_{\left( {2k + 1} \right)}^{\left( {2k + 2} \right)}  
= \left( {2k + 2} \right)^x  - \left( {2k + 1} \right)^x  = x\left( {2k + 1 + \lambda } \right)^{x - 1} 
\label{(B.10)}
\end{equation}
}
\begin{eqnarray}
\; & \; & 
\left| {\frac{{y\sin \left[ {y\log \left( {2m + 1 + \lambda } \right)} \right]}}{{\left( 
{2m + 1} \right)^x \left( {2m + 1 + \lambda } \right)}} + \frac{{x\cos \left[ {y\log \left( 
{2m + 2} \right)} \right]\left( {2m + 1 + \lambda } \right)^{x - 1} }}{{\left( {2m + 1} \right)^x 
\left( {2m + 2} \right)^x }}} \right|   
\nonumber \\
& \le & \left| {\frac{{y\sin \left[ {y\log \left( {2m + 1 + \lambda } \right)} \right]}}
{{\left( {2m + 1} \right)^x \left( {2m + 1 + \lambda } \right)}}} \right| + \left| 
{\frac{{x\cos \left[ {y\log \left( {2m + 2} \right)} \right]\left( {2m + 1 + \lambda } \right)^{x - 1} }}
{{\left( {2m + 1} \right)^x \left( {2m + 2} \right)^x }}} \right|   
\nonumber \\
& \le & \left| {\frac{y}{{\left( {2m + 1} \right)^{1 + x} }}} \right| 
+ \left| {\frac{x}{{\left( {2m + 1} \right)^{1 + x} }}} \right|.  
\label{(B.11)}
\end{eqnarray}
By substituting in the series and taking the limits for $k\to\infty$ we prove the thesis
for any $x>0$. An analogous proof holds for the $G_{\alpha}$ functions.

Hence we obtain
\vskip 0.3cm
\noindent
{\bf Corollary 2 (Row and column sum)}
\begin{eqnarray}
\; & \; &
r_{n_2}=F_{1}(x,y){(-1)^{n_{2}} \over n_{2}^{x}} \cos(y \log n_{2})
+G_{1}(x,y){(-1)^{n_{2}}\over n_{2}^{x}} \sin(y \log n_{2}), 
\nonumber \\
& \; & c_{n_{1}}=F_{2}(x,y){(-1)^{n_{1}}\over n_{1}^{x}} \cos(y \log n_{1})
\nonumber \\
&+& G_{2}(x,y){(-1)^{n_{1}}\over n_{1}^{x}} \sin(y \log n_{1}).   
\label{(B.12)}
\end{eqnarray}
\vskip 0.3cm
\noindent
{\it Proof}. We point out that
\begin{eqnarray}
\; & \; &
r_{n_2 }  = \frac{{\left( { - 1} \right)^{n_2 } }}{{n_2^x }}\sum\limits_{n_1  = 1}^\infty  
{\frac{{\left( { - 1} \right)^{n_1 } }}{{n_1^x }}\cos \left( {y\log \frac{{n_2 }}{{n_1 }}} \right)}   
\nonumber \\ 
&=& \frac{{\left( { - 1} \right)^{n_2 } }}{{n_2^x }}\cos \left( {y\log n_2 } \right)\underbrace 
{\sum\limits_{n_1  = 1}^\infty  {\frac{{\left( { - 1} \right)^{n_1 } }}{{n_1^x }}
\cos \left( {y\log n_1 } 
\right)} }_{F_1 \left( {x,y} \right)} 
\nonumber \\
&+& \frac{{\left( { - 1} \right)^{n_2 } }}{{n_2^x }}
\sin \left( {y\log n_2 } \right)\underbrace {\sum\limits_{n_1  = 1}^\infty  
{\frac{{\left( { - 1} \right)^{n_1 } }}{{n_1^x }}\sin \left( 
{y\log n_1 } \right)} }_{G_1 \left( {x,y} \right)},
\label{(B.13)} 
\end{eqnarray}
\begin{eqnarray}
\; & \; &
c_{n_1 }  = \frac{{\left( { - 1} \right)^{n_1 } }}{{n_1^x }}\sum\limits_{n_2  = 1}^\infty  
{\frac{{\left( { - 1} \right)^{n_2 } }}{{n_2^x }}\cos \left( {y\log \frac{{n_2 }}{{n_1 }}} \right)}   
\nonumber \\
&=& \frac{{\left( { - 1} \right)^{n_1 } }}{{n_1^x }}\cos \left( {y\log n_1 } \right)\underbrace 
{\sum\limits_{n_2  = 1}^\infty  {\frac{{\left( { - 1} \right)^{n_2 } }}{{n_2^x }}
\cos \left( {y\log n_2 } \right)} }_{F_2 \left( {x,y} \right)} 
\nonumber \\
&+& \frac{{\left( { - 1} 
\right)^{n_1 } }}{{n_1^x }}\sin \left( {y\log n_1 } \right)\underbrace 
{\sum\limits_{n_2  = 1}^\infty  
{\frac{{\left( { - 1} \right)^{n_2 } }}{{n_2^x }}\sin \left( {y\log n_2 } 
\right)} }_{G_2 \left( {x,y} \right)}. 
\label{(B.14)}
\end{eqnarray}
Now we can prove
\vskip 0.3cm
\noindent
{\bf Theorem 4 (Pringsheim convergence of an alternating double series)}. The
double series defined for any $x+iy \in {\mathbb S}^+$
\begin{equation}
\sum\limits_{n_2  \ne n_1 } {\frac{{\left( { - 1} \right)^{n_1  + n_2 } }}{{\left( {n_2 n_1 } 
\right)^x }}\cos \left( {y\log \frac{{n_2 }}{{n_1 }}} \right)} 
\label{(B.15)}
\end{equation}
is Pringsheim-convergent.
\vskip 0.3cm
\noindent
{\it Proof}. In light of (B.6) and Theorem 3, for any $\varepsilon>0$ we can find a 
$\nu_{\varepsilon}$ such that, for any $p,q>\nu_{\varepsilon}$, we have
\begin{equation}
\left | \sum_{n_{2}=1}^{q}a_{n_{1} \; n_{2}}-c_{n1} \right| 
= \left | s_{n_{1},q}-c_{n_{1}} \right | 
< \varepsilon, \;
\left | \sum_{n_{1}=1}^{p} a_{n_{1} \; n_{2}}-r_{n_{2}} \right |
= \left | s_{p,n_{2}}-r_{n_{2}} \right | < \varepsilon.
\label{(B.16)}
\end{equation}
At this stage, inspired by Appendix A, we recognize that by a suitable increase of the
limit of summation we can reduce in a convenient way the $\varepsilon$ upper bound; 
in particular, we can find two numbers $l=l(p) \in {\mathbb N}$ and
$m=m(q) \in {\mathbb N}$ such that
\begin{equation}
\left | \sum_{n_{2}=1}^{q+l}a_{n_{1} \; n_{2}}-c_{n_{1}} \right |
=\left | s_{n_{1},q+l}-c_{n_{1}} \right |
< {\varepsilon \over \zeta(2) p^{2}}
\leq {\varepsilon \over \zeta(2) n_{1}^{2}} \;
\forall n_{1} \in \left \{ 1,...,p \right \},
\label{(B.17)}
\end{equation}
\begin{equation}
\left | \sum_{n_{1}=1}^{p+m}a_{n_{1} \; n_{2}}-r_{n_{2}} \right |
=\left | s_{p+m,n_{2}}-r_{n_{2}} \right |
< {\varepsilon \over \zeta(2) q^{2}}
\leq {\varepsilon \over \zeta(2) n_{2}^{2}} \;
\forall n_{2} \in \left \{ 1,...,q \right \}.
\label{(B.18)}
\end{equation}
Hence we have, summing over\footnote{For rows and columns the Pringsheim convergence 
explicitly asks that $p,q>\nu_{\varepsilon}$ are independent of each other. This condition 
can be expressed by requiring that $\left|p-q\right|<\infty$.} $n_1$ and $n_2$
\begin{eqnarray}
\; & \; &
\left| {\sum\limits_{n_1  = 1}^{p + m} {\left( {s_{n_1 ,q + l}  - c_{n_1 } } \right)} } \right| 
\leq \sum\limits_{n_1  = 1}^{p + m} {\left| {s_{n_1 ,q + l}  - c_{n_1 } } \right|}  
< \sum\limits_{n_1  = 1}^{p + m} {\frac{\varepsilon }{{\zeta \left( 2 \right)n_1^2 }}}  
< \varepsilon ,   
\nonumber \\
& \; & \left| {\sum\limits_{n_2  = 1}^{q + l} {\left( {s_{p + l,n_2 }  - r_{n_2 } } \right)} } \right| 
\leq \sum\limits_{n_2  = 1}^{q + l} {\left| {s_{p + l,n_2 }  - r_{n_2 } } \right|}  
\nonumber \\
& < & \sum\limits_{n_2  = 1}^{q + l} {\frac{\varepsilon }{{\zeta \left( 2 \right)n_2^2 }}}  
< \varepsilon .  
\label{(B.19)}
\end{eqnarray}
By virtue of the previous observation, for any $\varepsilon>0$, we can find 
a $\mu_{\varepsilon}$ such that, for positive integers $u,v > \mu_{\varepsilon}$,
\begin{equation}
{\left| {\sum\limits_{n_1  = 1}^{p+m+u} {\frac{{\left( { - 1} \right)^{n_1 } }}{{n_1^x }}
\cos \left( {y\log n_1 } \right)}  - F_1 \left( {x,y} \right)} \right| < {\varepsilon \over 4} },
\label{(B.20)}
\end{equation}
\begin{equation}
{\left| {\sum\limits_{n_2  = 1}^{q+l+v} {\frac{{\left( { - 1} 
\right)^{n_2 } }}{{n_2^x }}\cos \left( {y\log n_2 } \right)}  - F_2 \left( {x,y} 
\right)} \right| < {\varepsilon \over 4} },  
\label{(B.21)}
\end{equation}
and analogously we can write
\begin{equation}
{\left| {\sum\limits_{n_1  = 1}^{p+m+u} {\frac{{\left( { - 1} \right)^{n_1 } }}{{n_1^x }}
\sin \left( {y\log n_1 } \right)}  - G_1 \left( {x,y} \right)} \right| < \frac{\varepsilon}{4} },
\label{(B.22)}
\end{equation}
\begin{equation}
{\left| {\sum\limits_{n_2  = 1}^{q+l+v} {\frac{{\left( { - 1} \right)^{n_2 } }}
{{n_2^x }}\sin \left( {y\log n_2 } \right)}  - G_2 \left( {x,y} \right)} \right| < \frac{\varepsilon}{4} }.
\label{(B.23)}  
\end{equation}
Thus, from the previous relations we find
\begin{eqnarray}
\; & \; &
\left| {\sum\limits_{n_2  = 1}^{q+l+v} {r_{n_2 } }  - \sum\limits_{n_1  = 1}^{p+m+u} {c_{n_1 } } } \right|  
\nonumber \\
& < & 
\left| F_{1}(x,y) \sum_{n_{2}=1}^{q+l+v}
{(-1)^{n_{2}}\over n_{2}^{x}} \cos(y \log n_{2}) \right .
\nonumber \\
&-& \left . F_{2}(x,y) \sum_{n_{1}=1}^{p+m+u}
{(-1)^{n_{1}}\over n_{1}^{x}}
\cos(y \log n_{1}) \right |
\nonumber \\
&+&\left| G_{1}(x,y) \sum_{n_{2}=1}^{q+l+v}
{(-1)^{n_{2}}\over n_{2}^{x}} \sin(y \log n_{2}) \right .
\nonumber \\
&-& \left . G_{2}(x,y) \sum_{n_{1}=1}^{p+m+u}
{(-1)^{n_{1}}\over n_{1}^{x}}
\sin(y \log n_{1}) \right |
\nonumber \\
& < & \left| {F_1 \left( {x,y} \right)\left( {F_2 \left( {x,y} \right) + \frac{\varepsilon}{4} } 
\right) - F_2 \left( {x,y} \right)\left( {F_1 \left( {x,y} \right) 
- \frac{\varepsilon}{4} } \right)} \right|  
\nonumber \\
&+& \left| {G_1 \left( {x,y} \right)\left( {G_2 \left( {x,y} \right) + \frac{\varepsilon}{4} } 
\right) - G_2 \left( {x,y} \right)\left( {G_1 \left( {x,y} \right) - \frac{\varepsilon}{4} } 
\right)} \right| 
\nonumber \\
& < & {\varepsilon \over 2}+{\varepsilon \over 2}
< \varepsilon.
\label{(B.24)}
\end{eqnarray}
At this stage, upon evaluating the difference between the sums $S_{a,b}$ and
${\tilde S}_{b,a}$ with the order of summation inverted, i.e.
 
\begin{eqnarray}
\; & \; &
\left| {S_{p+m+u,q+l+v}  - {\tilde S}_{q+l+v,p+m+u} } \right| 
= \left| {\sum\limits_{n_1  = 1}^{p+m+u} 
{s_{n_1 ,q+l+v} }  - \sum\limits_{n_2  = 1}^{q+l+v} {s_{p+m+u,n_2 } } } \right|   
\nonumber \\
& \leq & \left | \sum_{n_{1}=1}^{p+m+u}
\Bigr(s_{n_{1},q+l+v}-c_{n_{1}}\Bigr) \right |
+\left | \sum_{n_{2}=1}^{q+l+v} \Bigr(r_{n_{2}}-s_{p+m+u,n_{2}}\Bigr) \right |
\nonumber \\
&+& \left | \sum_{n_{1}=1}^{p+m+u}c_{n_{1}}
-\sum_{n_{2}=1}^{q+l+v} r_{n_{2}} \right | 
\nonumber \\
& \leq & \sum\limits_{n_1  = 1}^{p+m+u} {\left| {s_{n_1 ,q+l+v}  - c_{n_1 } } \right|}  
+ \sum\limits_{n_2  = 1}^{q+l+v} {\left| {s_{p+m+u,n_2 }  - r_{n_2 } } \right|}  
\nonumber \\
& + & \left| {\sum\limits_{n_1  = 1}^{p+m+u} {c_{n_1 } }  - \sum\limits_{n_2  = 1}^{q+l+v} 
{r_{n_2 } } } \right|
\nonumber \\ 
& < & 3\varepsilon ,  
\label{(B.25)} 
\end{eqnarray}
we have the thesis. These conclusions remain true when $u,v \to\infty$.
\vskip 0.3cm
\noindent
{\bf Corollary 3}. For any $x+iy \in {\mathbb S}^+$ the series
\begin{equation}
\sum\limits_{k = 1}^\infty  {\sum\limits_{n = 1}^\infty  {\frac{{\left( { - 1} \right)^{k + 1} }}
{{\left[ {n\left( {n + k} \right)} \right]^x }}
\cos \left[ {y\log \left( {1 + \frac{k}{n}} \right)} \right]} } 		
\label{(B.26)}
\end{equation}
is Pringsheim-convergent.
\vskip 0.3cm
\noindent
{\it Proof}. For any $n,k \in {\mathbb N}$, setting $n_1\equiv n$ and $n_2\equiv n+k$ in the series
\begin{equation}
\sum\limits_{n_1  \ne n_2 } {\frac{{\left( { - 1} \right)^{n_1  + n_2 } }}{{\left( {n_1 n_2 } 
\right)^x }}\cos \left( {y\log \frac{{n_1 }}{{n_2 }}} \right)} , 
\label{(B.27)}
\end{equation}
its Pringsheim-convergence follows from the previous theorem.

\section{Maximal extension of the fundamental identity at $s=1$}
\setcounter{equation}{0}
 
We check the fundamental identity (4.3) for the value $x=1$:
In light of\footnote{By studying the limit
\[
\begin{array}{*{20}c}
{\mathop {\lim }\limits_{x \to 0}} \frac{{a^x  - 1}}{x} = \log a & {} &  
\Rightarrow  & {} & {\mathop {\lim }\limits_{x \to 1}} \frac{{a^{1 - x}  - 1}}{{1 - x}} 
= \log a , \\
\end{array}
\]
and bearing in mind that $\zeta$ has residue 1 at $s=1$ we have
\[
\begin{array}{*{20}c}

{\mathop {\lim}\limits_{x \to 1}} \frac{{a^{1 - x}  - 1}}{{1 - x}} = \log a & {} & {} & {} & {}  \\
{} & {} &  \Rightarrow  & {} & {\mathop  {\lim}\limits_{x \to 1}} \left( {1 - 2^{1 - x} } 
\right)\zeta(x) =  - {\mathop {\lim }\limits_{x \to 1}} 
\left( {2^{1 - x}  - 1} \right)\zeta (x) =   \\
{\mathop {\lim}\limits_{x \to 1}} \zeta (x)(x - 1) = 1 & {} & {} & {} & {}  \\
{} & {} & {} & {} & { =  -  {\mathop {\lim}\limits_{x \to 1}} \frac{{2^{1 - x}  - 1}}{{x - 1}} 
\cdot {\mathop {\lim} \limits_{x \to 1}} \zeta (x)\left( {x - 1} \right) = \log 2},  \\
\end{array}
\]
hence we can evaluate
\[
\begin{array}{*{20}c}
{\mathop {\lim}\limits_{x \to 0}} \sum_{k=1}^{\infty}  {\frac{{( - 1)^{k + 1} }}
{{(k+1)^x }}}  = \frac{1}{2} & {} & {\rm and} & {} 
& {\mathop {\lim}\limits_{x \to 1}} 
\sum_{k=1}^{\infty}  {\frac{{(-1)^{k + 1} }}{{(k+1)^x }}}   \\
\end{array} = 1 - \log 2 .
\]
}
\begin{equation}
{\mathop {\lim}\limits_{s \to 1}} \left( {1 - 2^{1 - s} } \right)\left( {1 - 2^{1 - \bar s} } 
\right)\zeta(s)\zeta({\bar s}) = (\log 2)^{2}, \; 
{\mathop {\lim}\limits_{x \to 1}} \zeta ({2x}) = \frac{{\pi ^2 }}{6},  
\label{(C.1)}
\end{equation}
considering the triangular convergence of the right-hand side of the equation which assures 
us of the Pringsheim convergence, checking that it is column- and row-convergent 
we obtain that the result is independent of the order of summation (Pringsheim Theorem): 
\begin{eqnarray}
\; & \; &
 \sum_{n_1  \ne n_2 } {\frac{{\left( { - 1} \right)^{n_1  + n_2 } }}{{n_1 n_2 }}}  
= \sum_{{\cal T}_r } {\left( {\begin{array}{*{20}c}
    *  & { - \frac{1}{2}} & {\frac{1}{3}} & { - \frac{1}{4}} & {\frac{1}{5}} & { - \frac{1}{6}} &  \cdots   \\
   { - \frac{1}{2}} &  *  & { - \frac{1}{6}} & {\frac{1}{8}} & { - \frac{1}{{10}}} &  \ddots  & {}  \\
   {\frac{1}{3}} & { - \frac{1}{6}} &  *  & { - \frac{1}{{12}}} &  \ddots  & {} & {}  \\
   { - \frac{1}{4}} & {\frac{1}{8}} & { - \frac{1}{{12}}} &  \ddots  & {} & {} & {}  \\
   {\frac{1}{5}} & { - \frac{1}{{10}}} &  \ddots  & {} &  *  & {} & {}  \\
   { - \frac{1}{6}} &  \ddots  & {} & {} & {} &  \ddots  & {}  \\
    \vdots  & {} & {} & {} & {} & {} &  \ddots   \\
\end{array}} \right)}  
\nonumber  \\  
&=& 2\left\{ {0 - \frac{1}{2} + \frac{1}{3} - \frac{1}{4} - \frac{1}{6} + \frac{1}{5} 
+ \frac{1}{8} - \frac{1}{6} - \frac{1}{{10}} - \frac{1}{{12}} +  \cdots } \right\} 
\nonumber \\
& \longrightarrow & \left(\log 2\right)^2-\frac{\pi^2}{6} \approx  - 1.164481, 
\label{(C.2)}
\end{eqnarray}
which is the correct result.

\end{appendix}

\section*{acknowledgements}
P. D'Isanto dedicates his work to his parents.
G. Es\-po\-si\-to is gra\-te\-ful to the Di\-par\-ti\-men\-to di Fi\-si\-ca 
``Et\-to\-re Pan\-ci\-ni'' for hos\-pi\-ta\-li\-ty a\-nd sup\-po\-rt.
The authors are grateful to Professor Vincenzo Ferone for his remarks that led
to the creation of sub-section 6.1, to Professor Luigi Rosa for encouragement,
and to Professor Jacques G\'elinas for correspondence.
Last, but not least, encouragement and advice from Professor Klaus Kirsten 
have provided strong motivation for completing our work.

\end{document}